\documentclass[11pt]{article}
\pdfoutput=1

\usepackage{fullpage,xcolor}
\usepackage{mathabx}
\usepackage{graphicx}
\graphicspath{{./},{./images/},}

\usepackage{pdfpages,hyperref}
\hypersetup{colorlinks=true, linkcolor=blue,  anchorcolor=blue, citecolor=blue, filecolor=blue, menucolor=blue, pagecolor=blue, urlcolor=blue}

\newcommand{\bul}{$\sqbullet$}

\usepackage{fancyhdr}
\begin{document}

\thispagestyle{empty}

\ 
\vspace{4cm}
\begin{center}
{\Huge \bf
Proceedings of the second\\
``international Traveling Workshop\\
on Interactions between Sparse models and Technology''\\[2mm]
\textsc{iTWIST'14}\\[.5cm]
}
{\Large L'Arsenal, Namur, Belgium}\\[2mm]
{August 27-29, 2014.}\\[1cm]
\end{center}

\newpage
\null
\vfill
\begin{center}
\includegraphics[width=\textwidth]{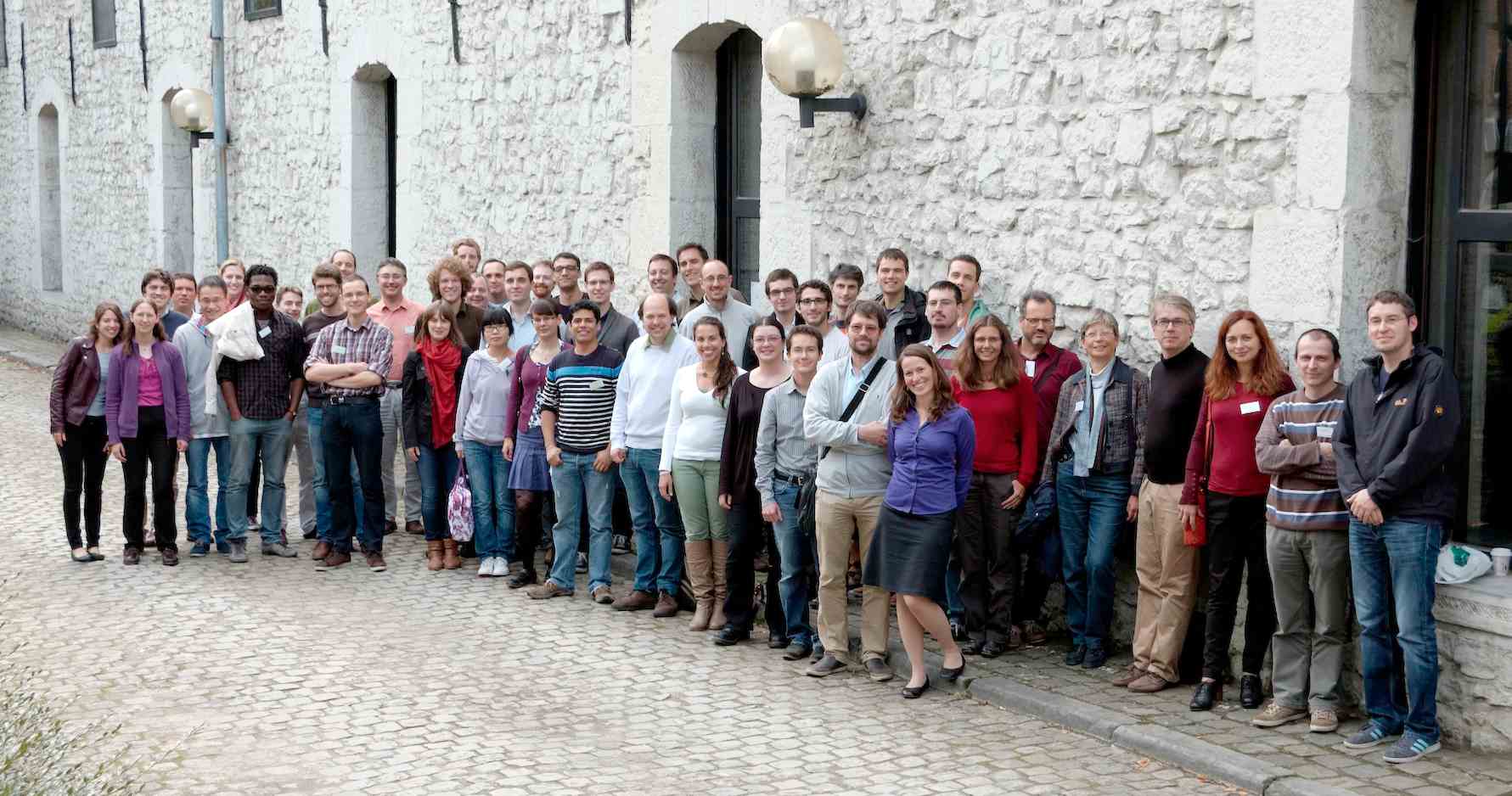}\\
iTWIST'14 Group Photo  
\end{center}
\vfill
\null

\newpage
\setcounter{page}{1}
\renewcommand{\thepage}{\roman{page}}

\begin{center}
\Large i\textsc{TWIST'14 Presentation}  
\end{center}

The second edition of the
``\emph{international - Traveling Workshop on Interactions between Sparse
models and Technology}'' (iTWIST) took place in the medieval and picturesque
town of Namur in Belgium. The workshop was conveniently located in ``The
Arsenal'' building within walking distance of both hotels and town center.
One implicit objective of this biennial workshop is to foster collaboration
between international scientific teams by disseminating ideas through
both specific oral/poster presentations and free discussions.
 
For this second edition, iTWIST'14 has gathered about 70 international
participants and has featured 9 invited talks, 10 oral presentations, and 14 posters on the following
themes, all related to the theory, application and generalization of
the ``sparsity paradigm'':\\[-10mm]
\begin{center}
\begin{tabular}{ll}
\bul~Sparsity-driven data sensing and processing&\bul~Union of low
dimensional subspaces\\
\bul~Beyond linear and convex inverse
problem&\bul~Matrix/manifold/graph sensing/processing\\
\bul~Blind inverse problems and dictionary learning&\bul~Sparsity and computational
neuroscience\\
\bul~Information theory, geometry and
randomness&\bul~Complexity/accuracy tradeoffs in numerical methods\\
\bul~Sparsity? What's next?&\bul~Sparse machine learning and inference
\end{tabular}
\end{center}

\begin{center}
\includegraphics[width=0.3\textwidth]{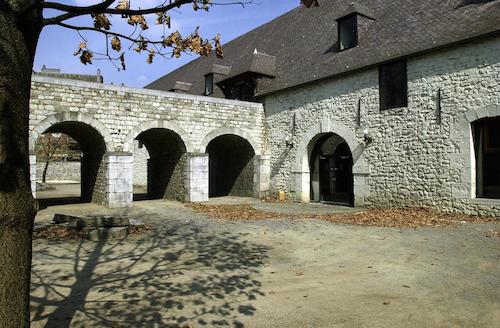}\hspace{.5mm}  
\includegraphics[width=0.3\textwidth]{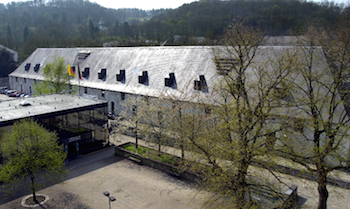}\hspace{.5mm}  
\includegraphics[width=0.3\textwidth]{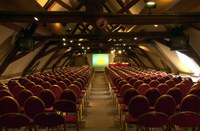}\\
\small Pictures of ``L'Arsenal'' Rue de Bruno, 11 5000 Namur - Belgique  
\end{center}

\textsc{Scientific Organizing Committee:}\\[4mm]
 \begin{tabular}{p{.5\textwidth}p{.5\textwidth}}
\bul~Laurent Jacques (UCL, Belgium)&\bul~Christophe De
Vleeschouwer (UCL, Belgium)\\
\bul~Yannick Boursier (CPPM, Aix-Marseille
U., France)&\bul~Prasad Sudhakar (UCL, Belgium)\\
\bul~Christine De Mol (ULB, Belgium)&\bul~Aleksandra Pizurica (Ghent U.,
Belgium)\\
\bul~Sandrine Anthoine (I2M, Aix-Marseille U., France)&\bul~Pierre
Vandergheynst (EPFL, Switzerland)\\
\bul~Pascal Frossard
(EPFL, Switzerland).\\  
\end{tabular}\\

\textsc{Local Organizing Committee:}\\[4mm]
\begin{tabular}{p{.5\textwidth}p{.5\textwidth}}
\bul~Laurent Jacques (UCL, Belgium)&
\bul~Christophe De Vleeschouwer (UCL, Belgium)\\
\bul~Jean Deschuyter (UCL, Belgium)&
\bul~Adriana Gonzalez (UCL, Belgium)\\
\bul~K\'evin Degraux (UCL, Belgium)&
\bul~St\'ephanie Gu\'erit (UCL, Belgium)\\
\bul~Pascaline Parisot (UCL, Belgium)&
\bul~Nicolas Matz (INSA, Toulouse, France)\\
\bul~Muhammad Arsalan (UCL, Belgium)&
\bul~C\'edric Verleysen (UCL, Belgium)\\
\bul~Amit Kumar KC (UCL, Belgium)&\\[4mm]
\end{tabular}

\noindent\textbf{Sponsors:} The iTWIST'14 organizing committee thanks the following sponsors for
their help and fundings.
\begin{center}
\includegraphics[width=\textwidth]{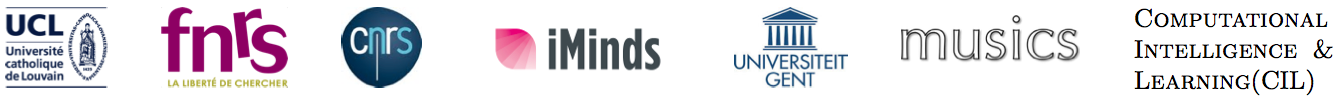}
\end{center}
\vfill
\null
\newpage
\section*{Table of contents}
\begin{itemize}
\footnotesize
\item[{\scriptsize (p.~\pageref*{pdf:Bilen})}] \hyperlink{Bilen.1}{C. Bilen (INRIA Rennes, France), Srdan Kitic (INRIA Rennes, France),
N. Bertin (IRISA CNRS, France), R. Gribonval (INRIA Rennes, France),\newline
``Sparse Acoustic Source Localization with Blind Calibration for Unknown Medium Characteristics''.}
\item[{\scriptsize  (p.~\pageref*{pdf:Boumal})}] \hyperlink{Boumal.1}{N. Boumal (UCL, Belgium), B. Mishra (ULg, Belgium), P.-A. Absil (UCL, Belgium), R. Sepulchre (Cambridge U., UK),\newline
``Manopt: a Matlab toolbox for optimization on manifolds''.}

\item[{\scriptsize  (p.~\pageref*{pdf:Bundervoet})}] \hyperlink{Bundervoet.1}{S. Bundervoet (VUB-ETRO,Belgium), C. Schretter (VUB-ETRO,Belgium),
A.Dooms (VUB-ETRO,Belgium), P. Schelkens (VUB-ETRO, Belgium),\newline
``Bayesian Estimation of Sparse Smooth Speckle Shape Models for Motion Tracking in Medical Ultrasound''.}

\item[{\scriptsize  (p.~\pageref*{pdf:Chabiron})}] \hyperlink{Chabiron.1}{O. Chabiron (U.Toulouse,France), F.Malgouyres (U.Toulouse,France),
J.-Y.,Tourneret (U.Toulouse,France), N. Dobigeon (U. Toulouse, France),\newline
``Learning a fast transform with a dictionary''.}

\item[{\scriptsize  (p.~\pageref*{pdf:Chainais})}] \hyperlink{Chainais.1}{P. Chainais (LAGIS / INRIA Lille, France), C. Richard (Lab. Lagrange, Nice, France),\newline
``A diffusion strategy for distributed dictionary learning''.}

\item[{\scriptsize  (p.~\pageref*{pdf:Cornelis})}] \hyperlink{Cornelis.1}{B. Cornelis (VUB-ETRO, Belgium), A. Dooms (VUB-ETRO, Belgium),
I. Daubechies (Duke U., USA), D. Dunson (Duke U., USA),\newline
``Bayesian crack detection in high resolution data''.}

\item[{\scriptsize  (p.~\pageref*{pdf:Dankova})}] \hyperlink{Dankova.1}{M. Dankova (Brno University of Technology) and P. Rajmic (Brno University of Technology, Czech Republic) \newline
``Compressed sensing of perfusion MRI''.}

\item[{\scriptsize  (p.~\pageref*{pdf:Degraux})}] \hyperlink{Degraux.1}{K. Degraux (UCL, Belgium), V. Cambareri (U.Bologna, Italy), B.Geelen (IMEC, Belgium),
L. Jacques (UCL, Belgium), G. Lafruit (IMEC, Belgium), G. Setti (U. Bologna, Italy),\newline
``Compressive Hyperspectral Imaging by Out-of-Focus Modulations and Fabry-P\'erot Spectral Filters''.}

\item[{\scriptsize  (p.~\pageref*{pdf:Determe})}] \hyperlink{Determe.1}{J.-F. Determe (ULB, Belgium), J. Louveaux
    (UCL, Belgium), F. Horlin (ULB, Belgium),\newline
``Filtered Orthogonal Matching Pursuit: Applications''.}

\item[{\scriptsize  (p.~\pageref*{pdf:Dremeau})}] \hyperlink{Dremeau.1}{A. Dr\'emeau (ESPCI ParisTech, France), P. H\'eas (INRIA, Rennes, France), C. Herzet (INRIA, Rennes, France),\newline
``Combining sparsity and dynamics: an efficient way''.}

\item[{\scriptsize  (p.~\pageref*{pdf:Duval})}] \hyperlink{Duval.1}{V. Duval (U. Paris-Dauphine), G. Peyr\'e (U. Paris-Dauphine),\newline
``Discrete vs. Continuous Sparse Regularization''.}

\item[{\scriptsize  (p.~\pageref*{pdf:Fawzi})}] \hyperlink{Fawzi.1}{A. Fawzi (LTS4, EPFL, Switzerland), M. Davies (U. Edinburgh, UK), P. Frossard (LTS4, EPFL, Switzerland),\newline
``Dictionary learning for efficient classification based on soft-thresholding''.}

\item[{\scriptsize  (p.~\pageref*{pdf:Gillis})}] \hyperlink{Gillis.1}{N. Gillis (U. Mons, Belgium), S. A. Vavasis (U. Waterloo, UK),\newline
``Semidefinite Programming Based Preconditioning for More Robust Near-Separable Nonnegative Matrix Factorization''.}

\item[{\scriptsize  (p.~\pageref*{pdf:Herzet})}] \hyperlink{Herzet.1}{C. Herzet (INRIA Rennes, France), C. Soussen (U. de Lorraine, Nancy, France),\newline
``Enhanced Recovery Conditions for OMP/OLS by Exploiting both Coherence and Decay''.}

\item[{\scriptsize  (p.~\pageref*{pdf:Kitic})}] \hyperlink{Kitic.1}{S. Kitic, (Inria, France), N. Bertin, (IRISA, France), R. Gribonval (Inria, France),\newline
``Wideband Audio Declipping by Cosparse Hard Thresholding''.}

\item[{\scriptsize  (p.~\pageref*{pdf:LeMagoarou})}] \hyperlink{LeMagoarou.1}{L. Le Magoarou (INRIA, Rennes, France), R. Gribonval (INRIA, Rennes, France),\newline
``Strategies to learn computationally efficient and compact dictionaries''.}

\item[{\scriptsize  (p.~\pageref*{pdf:Liang})}] \hyperlink{Liang.1}{J. Liang (GREYC, ENSICAEN, France), J. Fadili (GREYC, ENSICAEN, France), G. Peyr\'e (U. Paris-Dauphine, France),\newline
``Iteration-Complexity for Inexact Proximal Splitting Algorithms''.}

\item[{\scriptsize  (p.~\pageref*{pdf:Liutkus})}] \hyperlink{Liutkus.1}{A. Liutkus (Institut Langevin, France), D. Martina (Institut Langevin, France),
S. Gigan (Institut Langevin, France), L. Daudet (Institut Langevin, France),\newline
``Calibration and imaging through scattering media''.}

\item[{\scriptsize  (p.~\pageref*{pdf:Maggioni})}] \hyperlink{Maggioni.1}{M. Maggioni (Duke U.), S. Minsker (Duke U.), N. Strawn (Duke U.),\newline
``Multiscale Dictionary and Manifold Learning: Non-Asymptotic Bounds for the Geometric Multi-Resolution Analysis''.}

\item[{\scriptsize  (p.~\pageref*{pdf:Mory})}] \hyperlink{Mory.1}{C. Mory (iMagX Belgium and CREATIS U. Lyon France) and L. Jacques (UCLouvain, Belgium) \newline
``An application of the Chambolle-Pock algorithm to 3D + time tomography''.}

\item[{\scriptsize  (p.~\pageref*{pdf:Ngole})}] \hyperlink{Ngole.1}{F. Ngole (CEA-CNRS-Paris 7, France), Starck, Jean-Luc (CEA-CNRS-Paris 7, France),\newline
``Super-resolution method using sparse regularization for point spread function recovery''.}

\item[{\scriptsize  (p.~\pageref*{pdf:Schretter})}] \hyperlink{Schretter.1}{C. Schretter (VUB-ETRO, Belgium), I. Loris (ULB, Belgium),
A. Doom (VUB-ETRO, Belgium), P. Schelkens (VUB- ETRO, Belgium),\newline
``Total Variation Reconstruction From Quasi-Random Samples''.}

\item[{\scriptsize  (p.~\pageref*{pdf:Vaiter})}] \hyperlink{Vaiter.1}{S. Vaiter (U. Paris-Dauphine, France), M. Golbabaee (U. Paris-Dauphine, France),
J. Fadili (GREYC, ENSICAEN, France), G. Peyr\'e (U. Paris-Dauphine, France),\newline
``Model Selection with Piecewise Regular Gauges''.}

\item[{\scriptsize  (p.~\pageref*{pdf:Vukobratovic})}] \hyperlink{Vukobratovic.1}{D. Vukobratovic (U. of Novi Sad, Serbia), A. Pizurica (Ghent University, Belgium),\newline
``Adaptive Compressed Sensing Using Sparse Measurement Matrices''.}
\end{itemize}
\newpage
\setcounter{page}{1}
\renewcommand{\thepage}{\arabic{page}}
\label{pdf:Bilen}
\includepdf[offset=1cm -4mm,pages=-,link=true,linkname=Bilen,pagecommand={}]{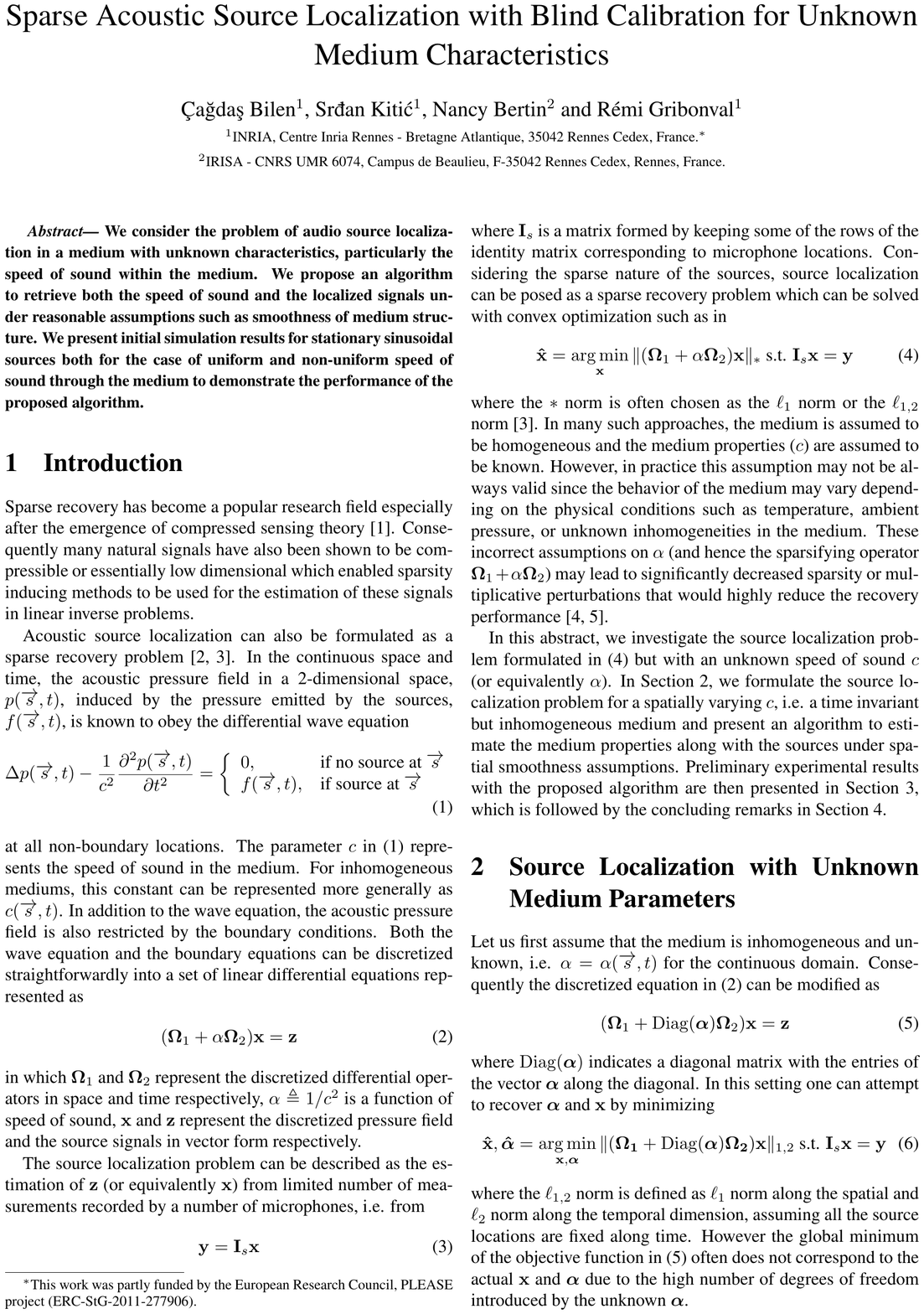}
\label{pdf:Boumal}
\includepdf[offset=1cm -4mm,pages=-,link=true,linkname=Boumal,pagecommand={}]{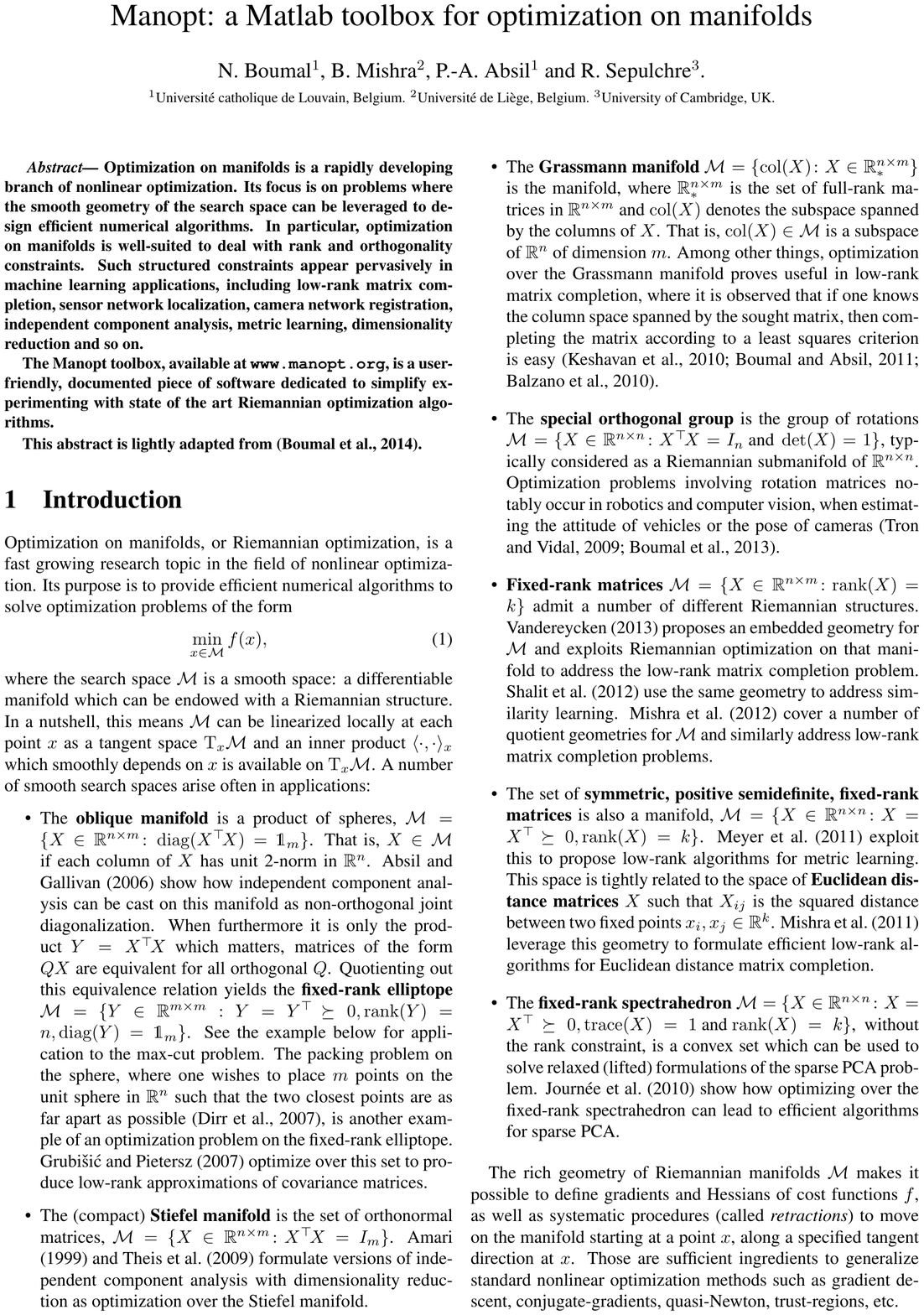}
\label{pdf:Bundervoet}
\includepdf[offset=1cm -4mm,pages=-,link=true,linkname=Bundervoet,pagecommand={}]{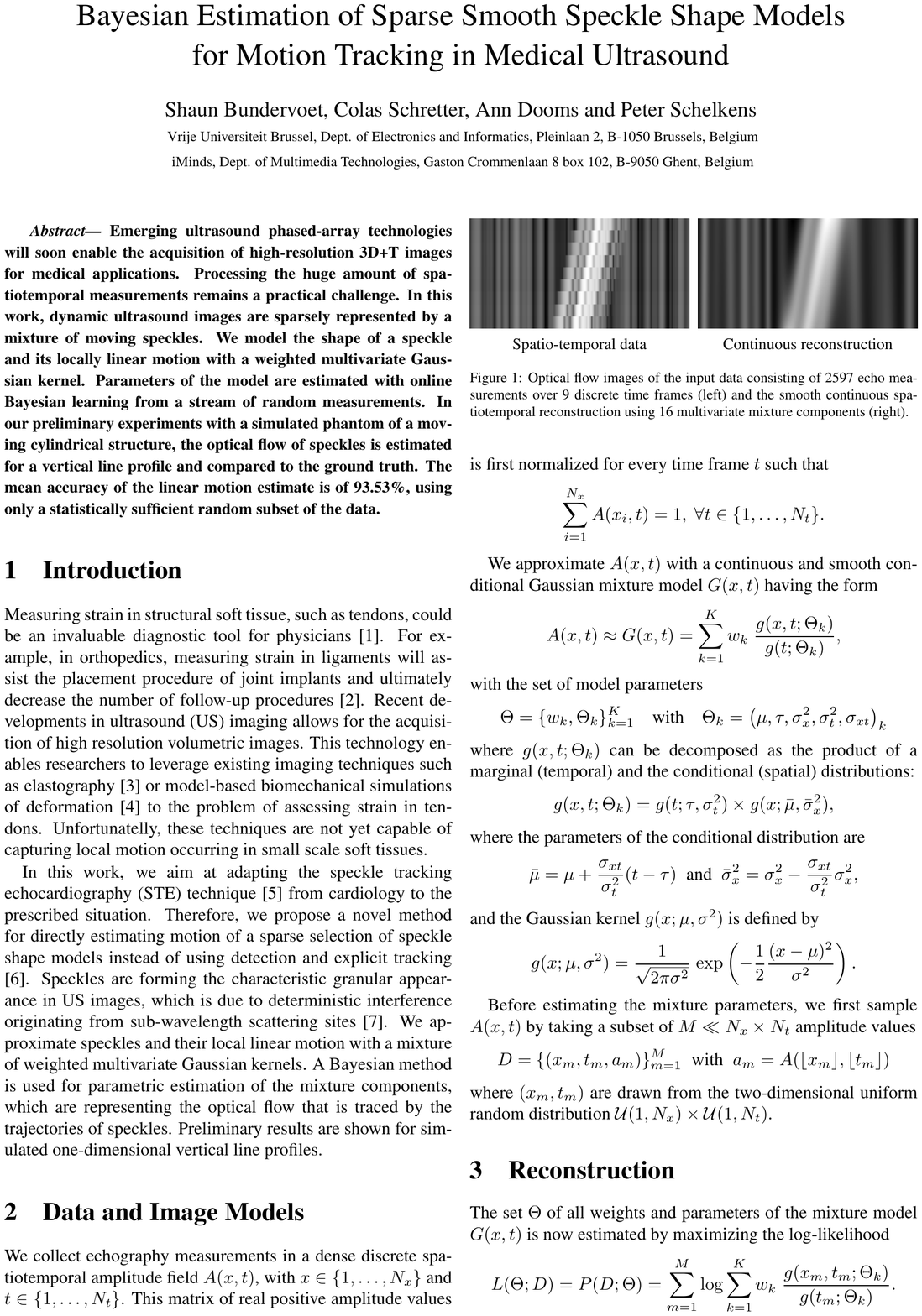}
\label{pdf:Chabiron}
\includepdf[offset=1cm -4mm,pages=-,link=true,linkname=Chabiron,pagecommand={}]{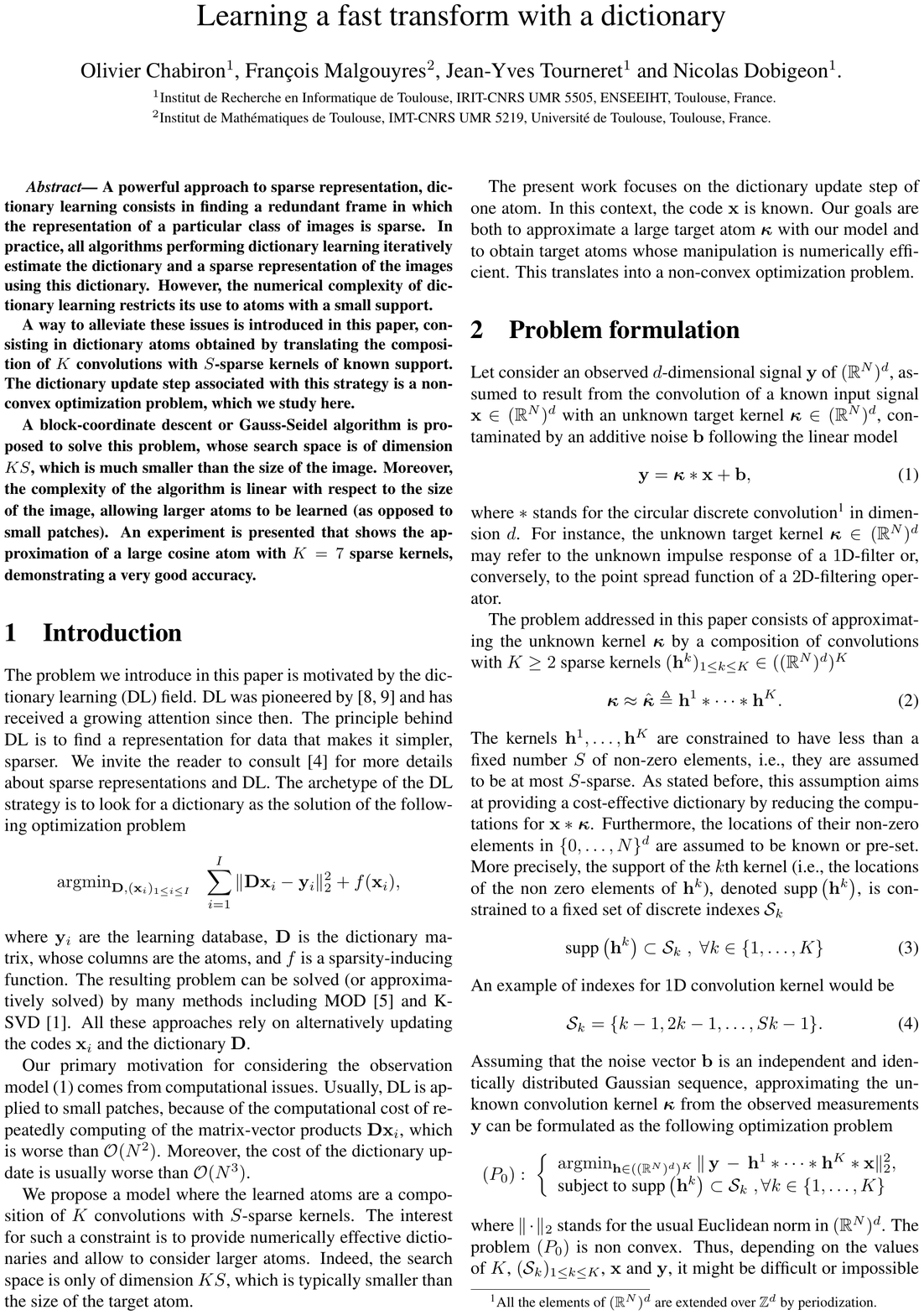}
\label{pdf:Chainais}
\includepdf[offset=1cm -4mm,pages=-,link=true,linkname=Chainais,pagecommand={}]{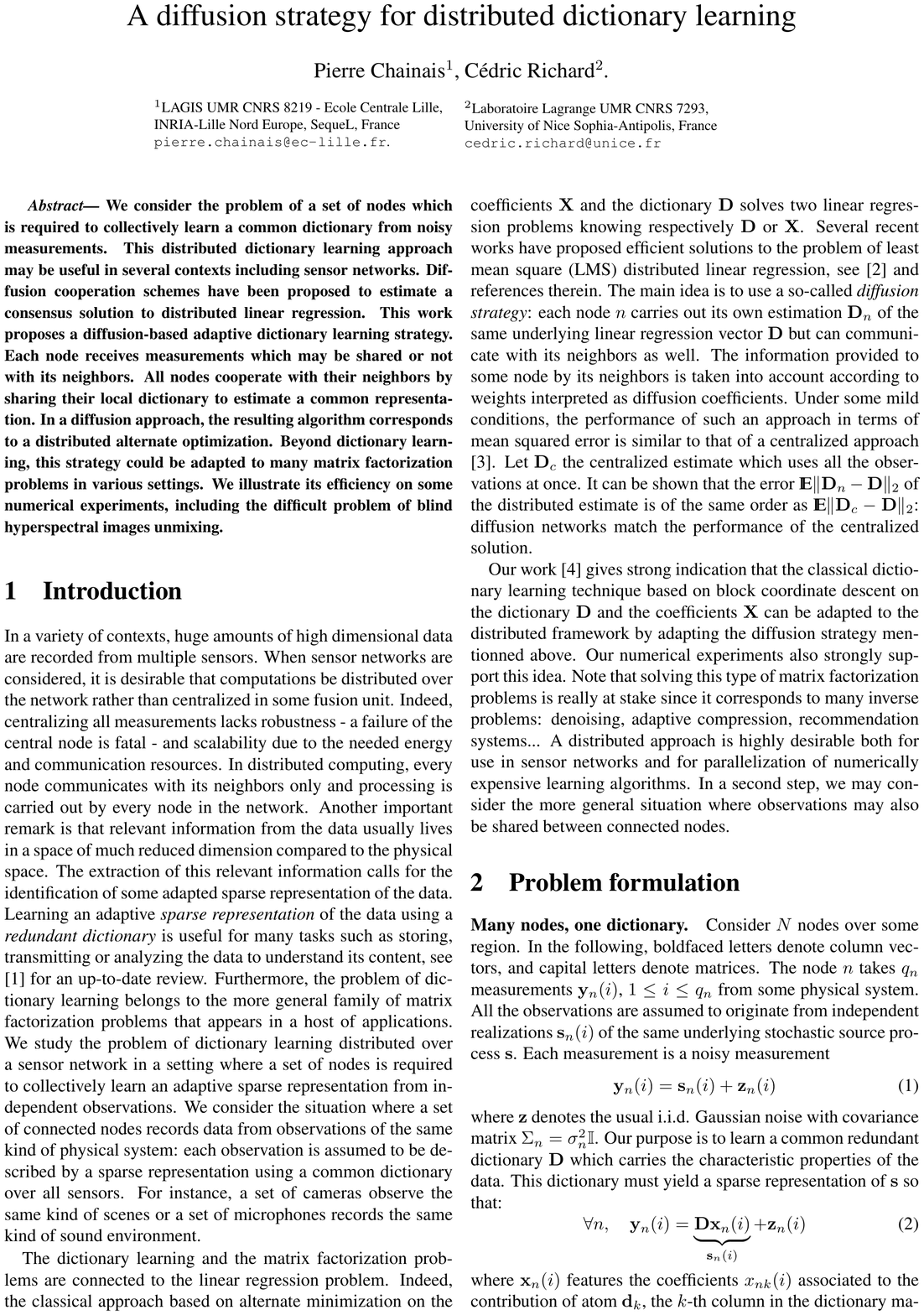}
\label{pdf:Cornelis}
\includepdf[offset=1cm -4mm,pages=-,link=true,linkname=Cornelis,pagecommand={}]{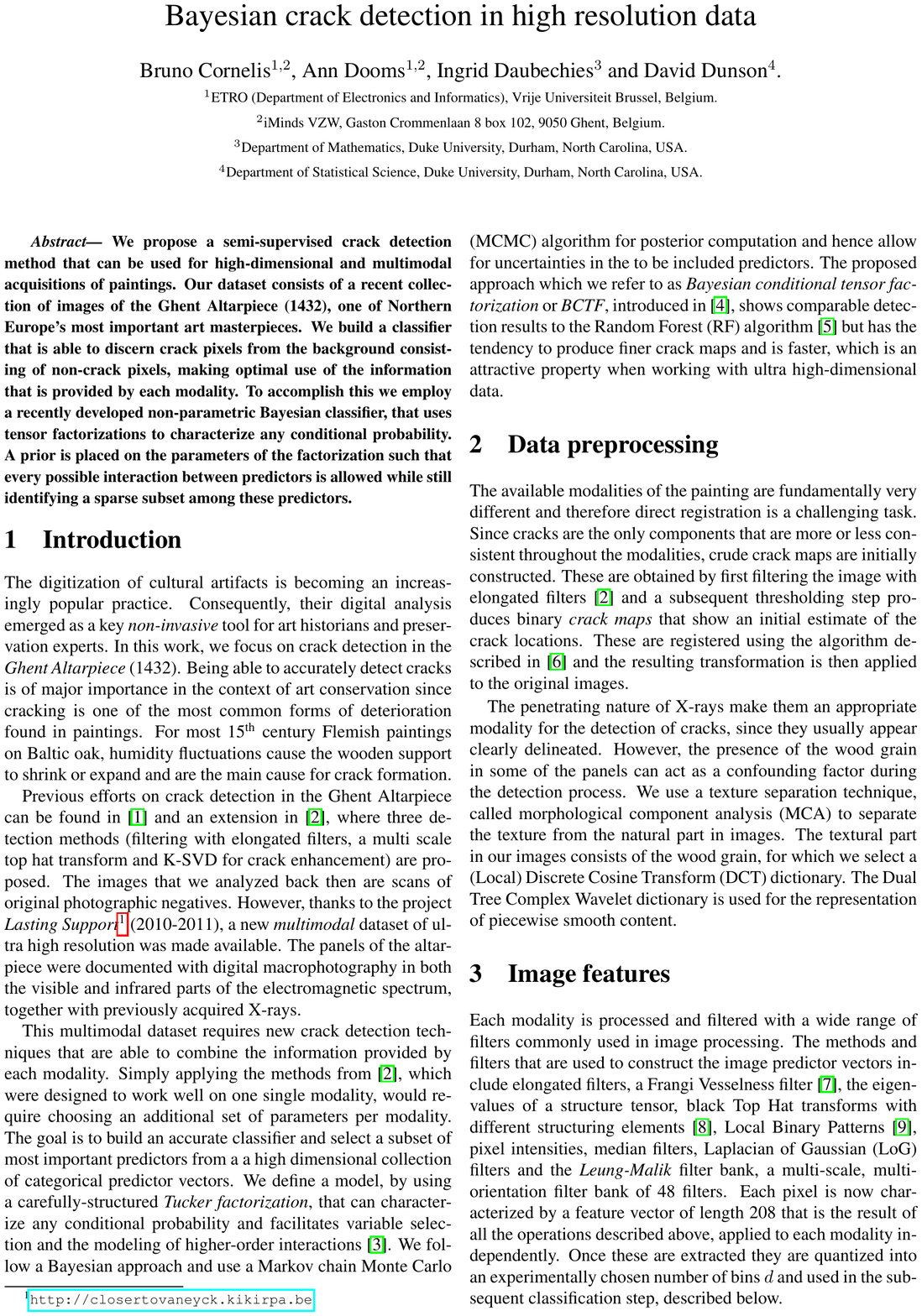}
\label{pdf:Dankova}
\includepdf[offset=1cm -4mm,pages=-,link=true,linkname=Dankova,pagecommand={}]{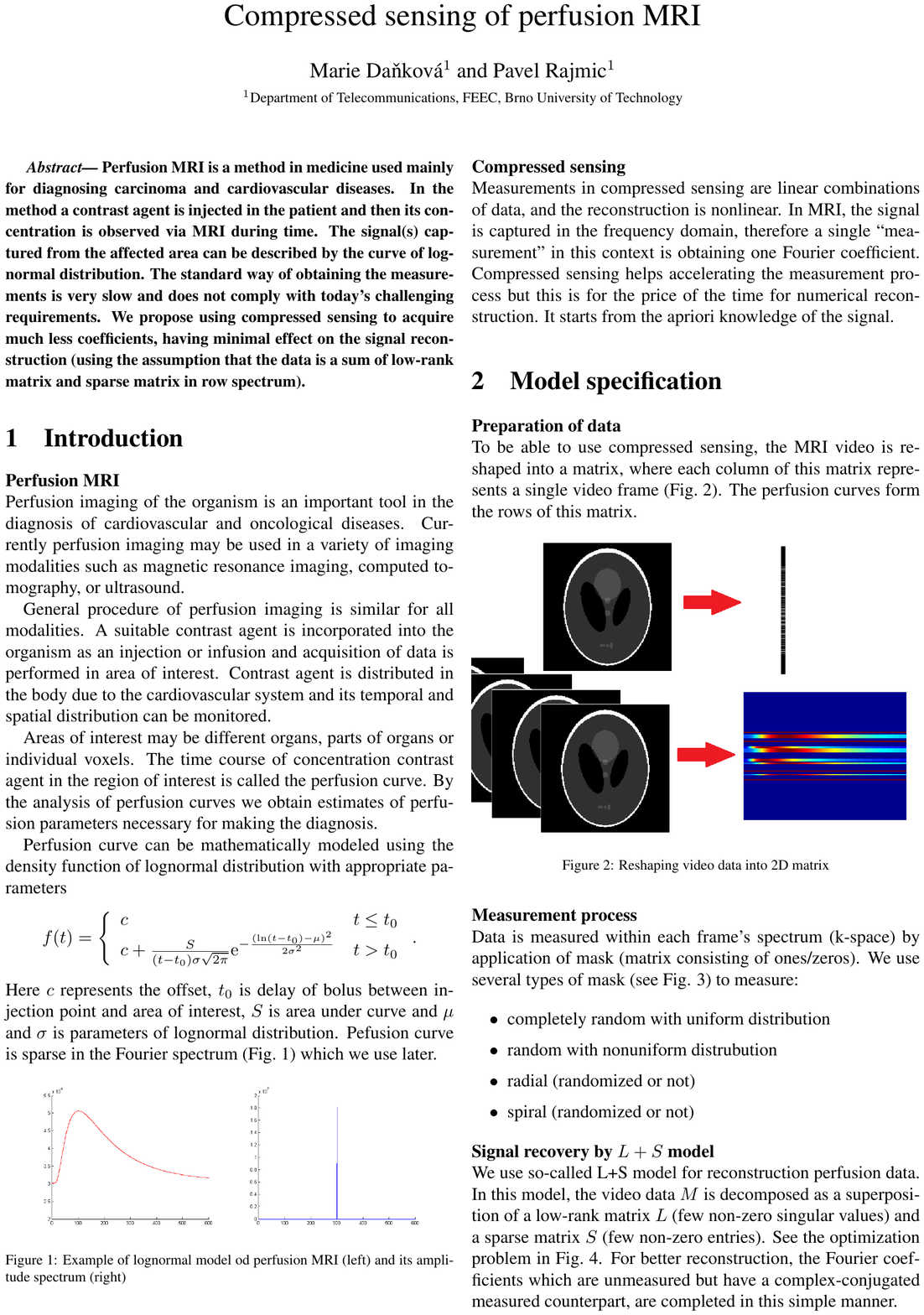}
\label{pdf:Degraux}
\includepdf[offset=1cm -4mm,pages=-,link=true,linkname=Degraux,pagecommand={}]{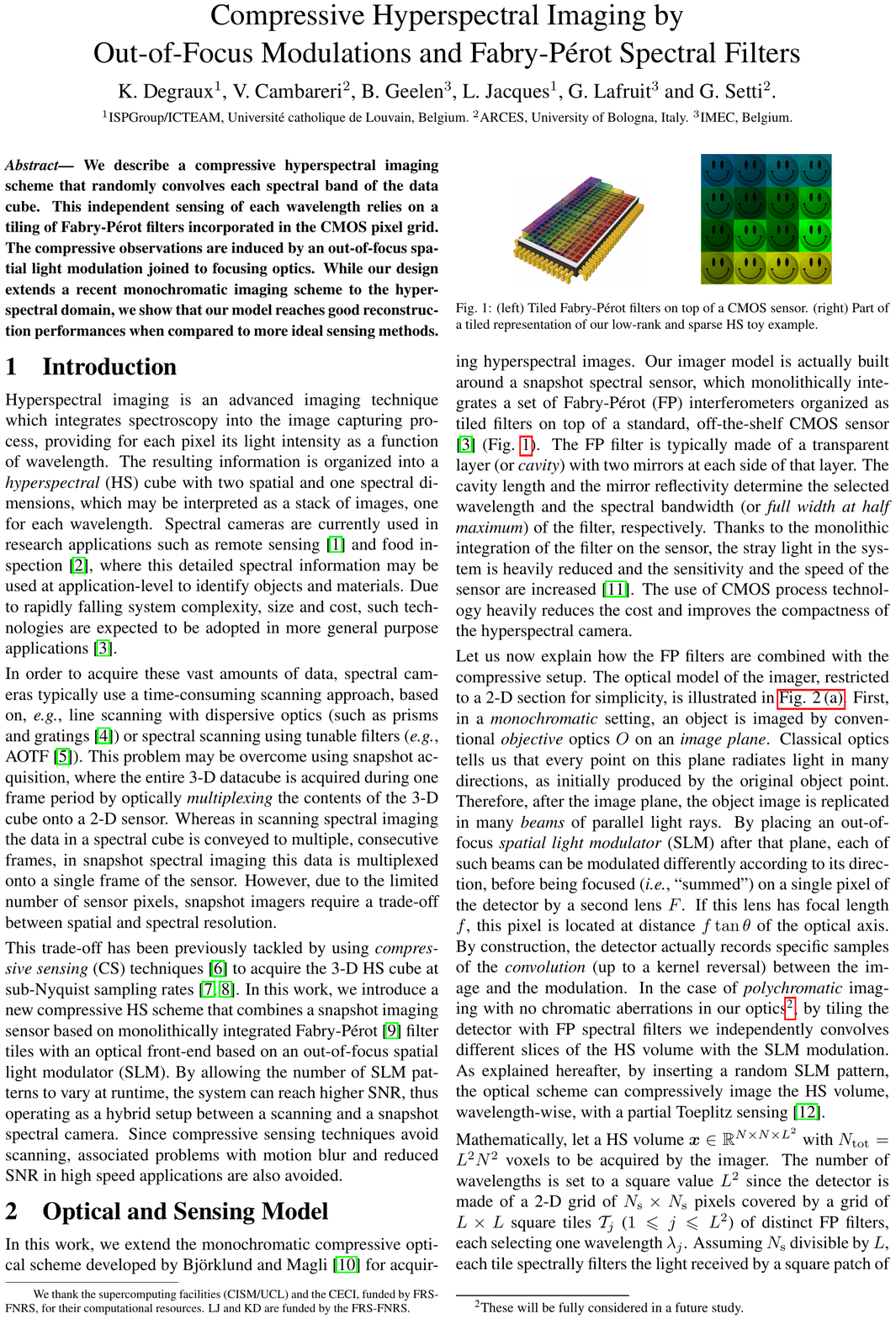}
\label{pdf:Determe}
\includepdf[offset=1cm -4mm,pages=-,link=true,linkname=Determe,pagecommand={}]{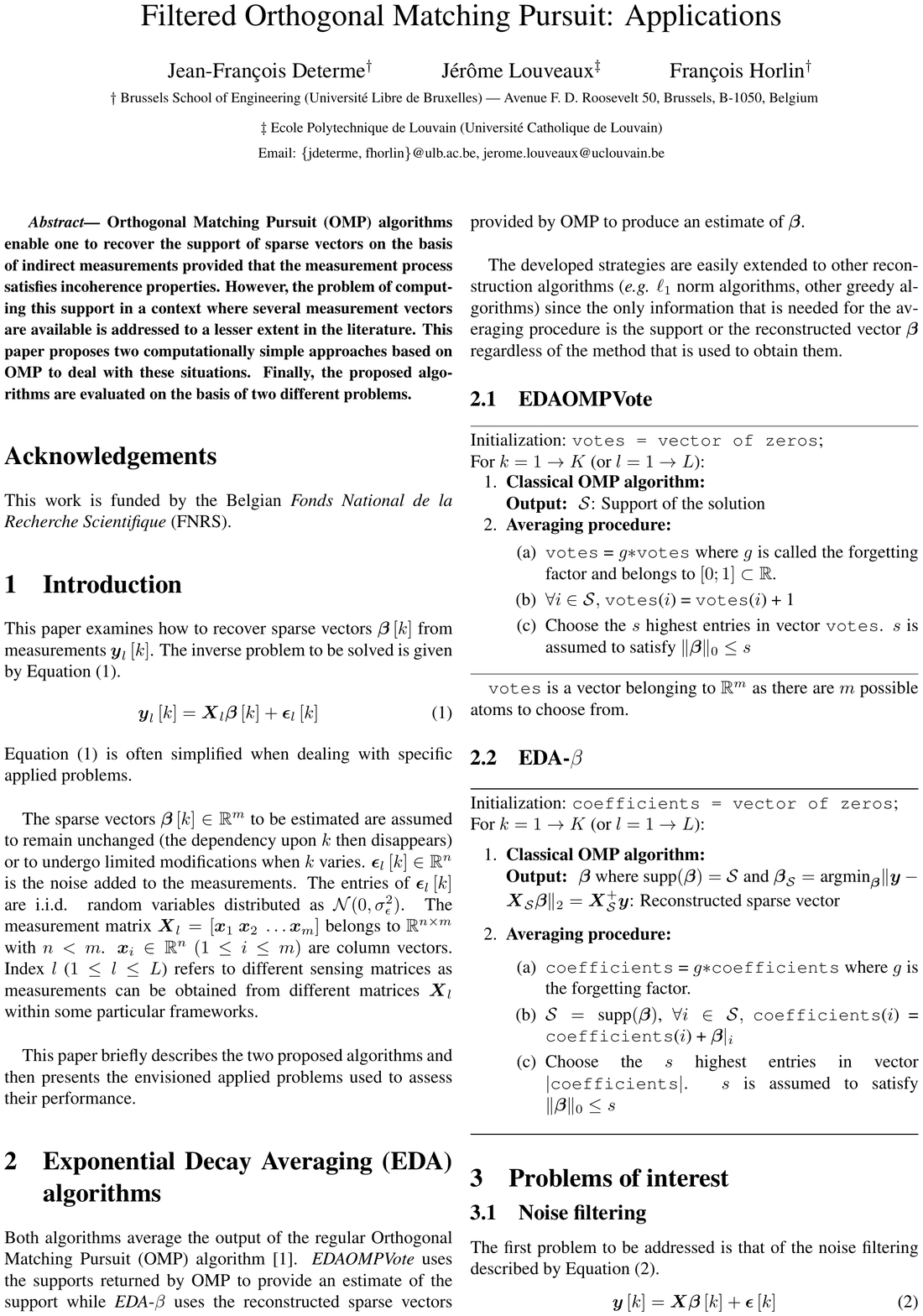}
\label{pdf:Dremeau}
\includepdf[offset=1cm -4mm,pages=-,link=true,linkname=Dremeau,pagecommand={}]{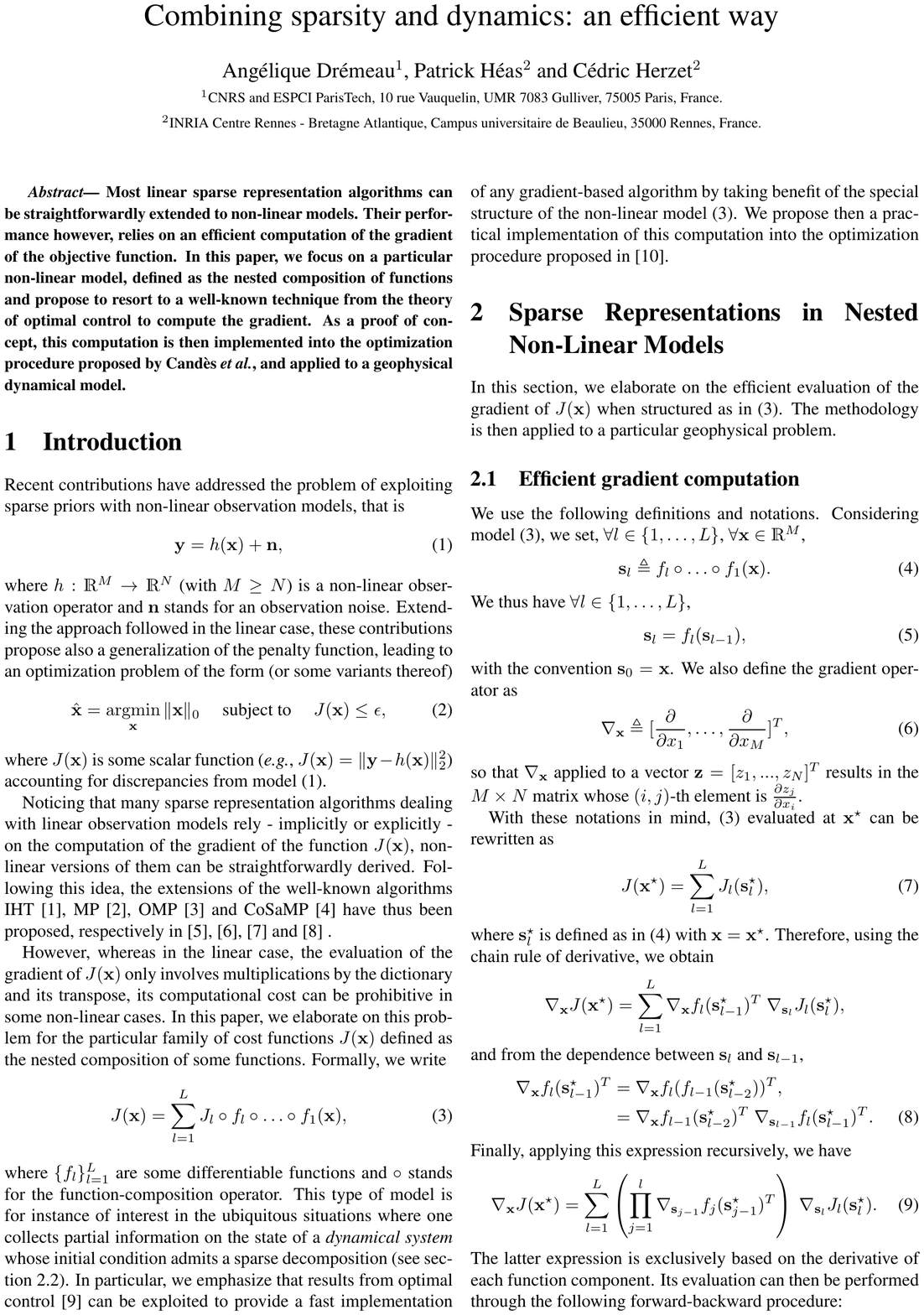}
\label{pdf:Duval}
\includepdf[offset=1cm -4mm,pages=-,link=true,linkname=Duval,pagecommand={}]{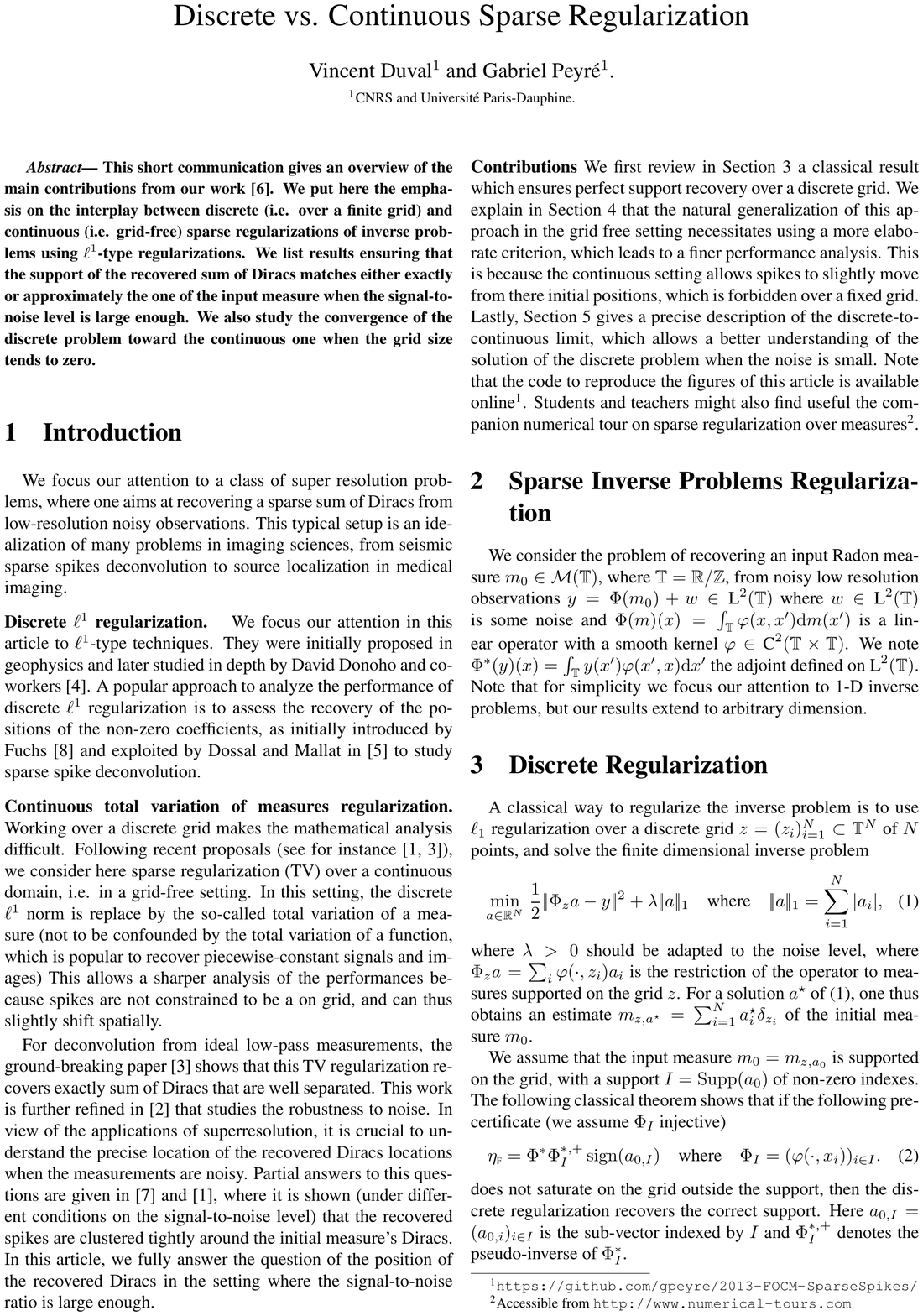}
\label{pdf:Fawzi}
\includepdf[offset=1cm -4mm,pages=-,link=true,linkname=Fawzi,pagecommand={}]{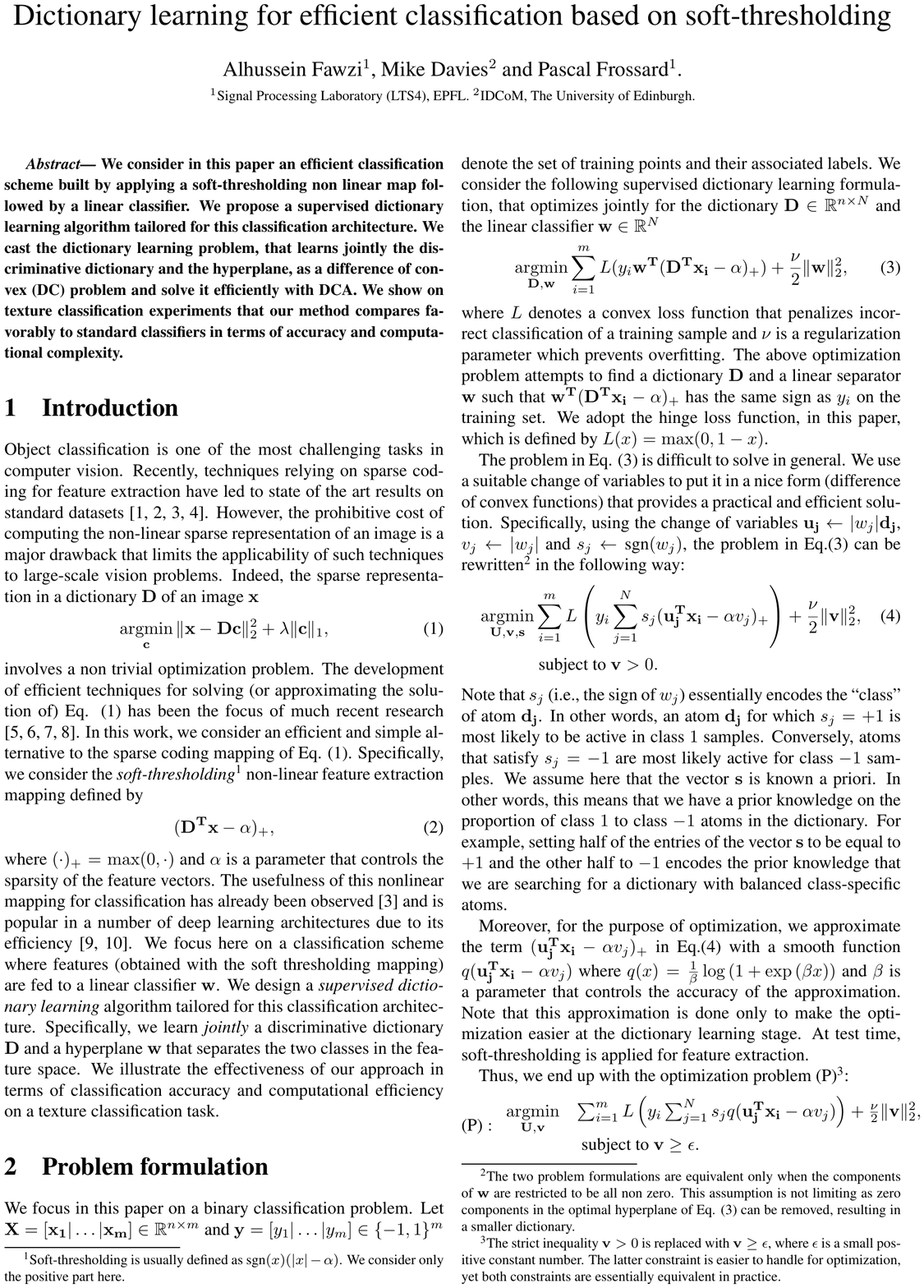}
\label{pdf:Gillis}
\includepdf[offset=1cm -4mm,pages=-,link=true,linkname=Gillis,pagecommand={}]{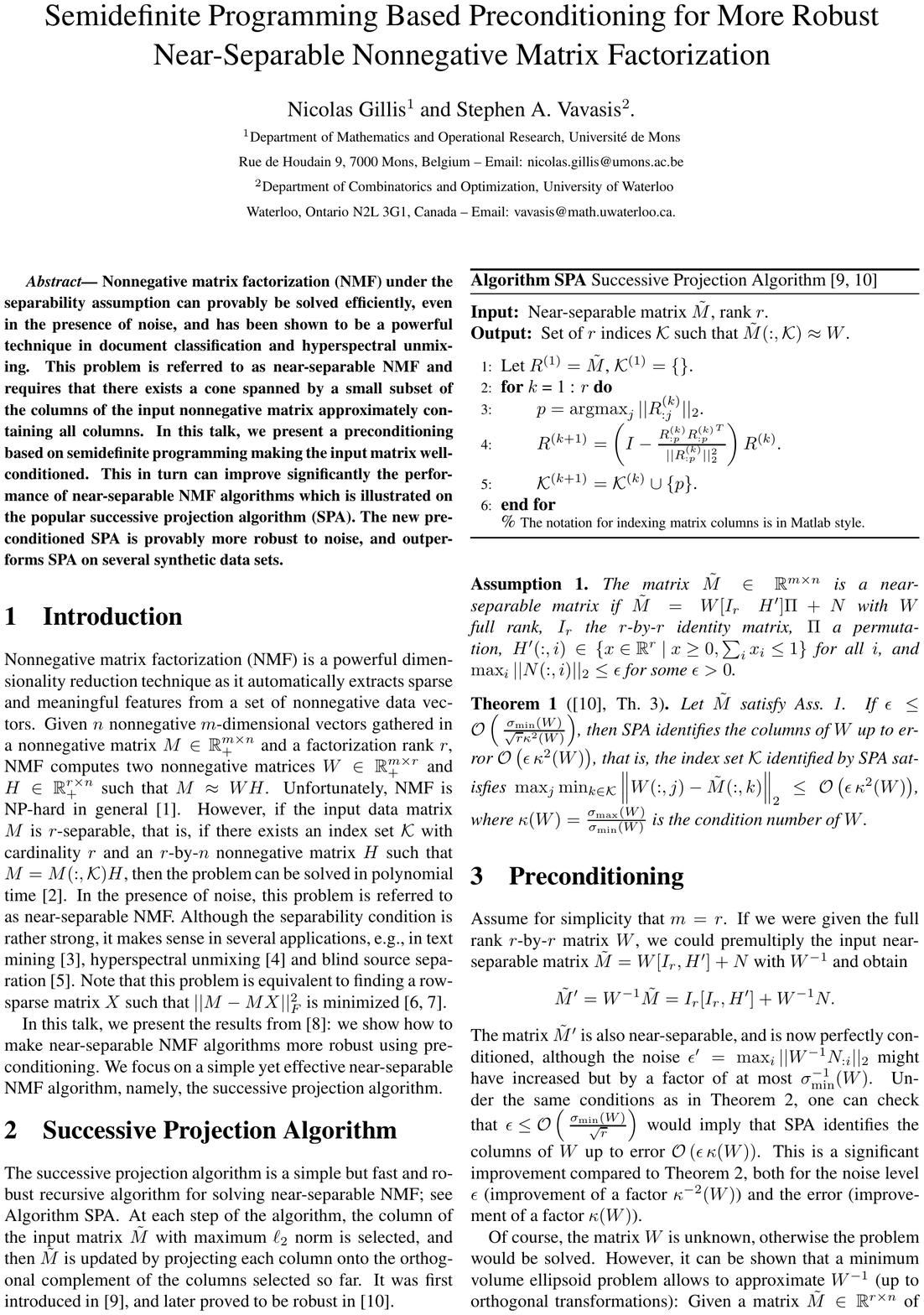}
\label{pdf:Herzet}
\includepdf[offset=1cm -4mm,pages=-,link=true,linkname=Herzet,pagecommand={}]{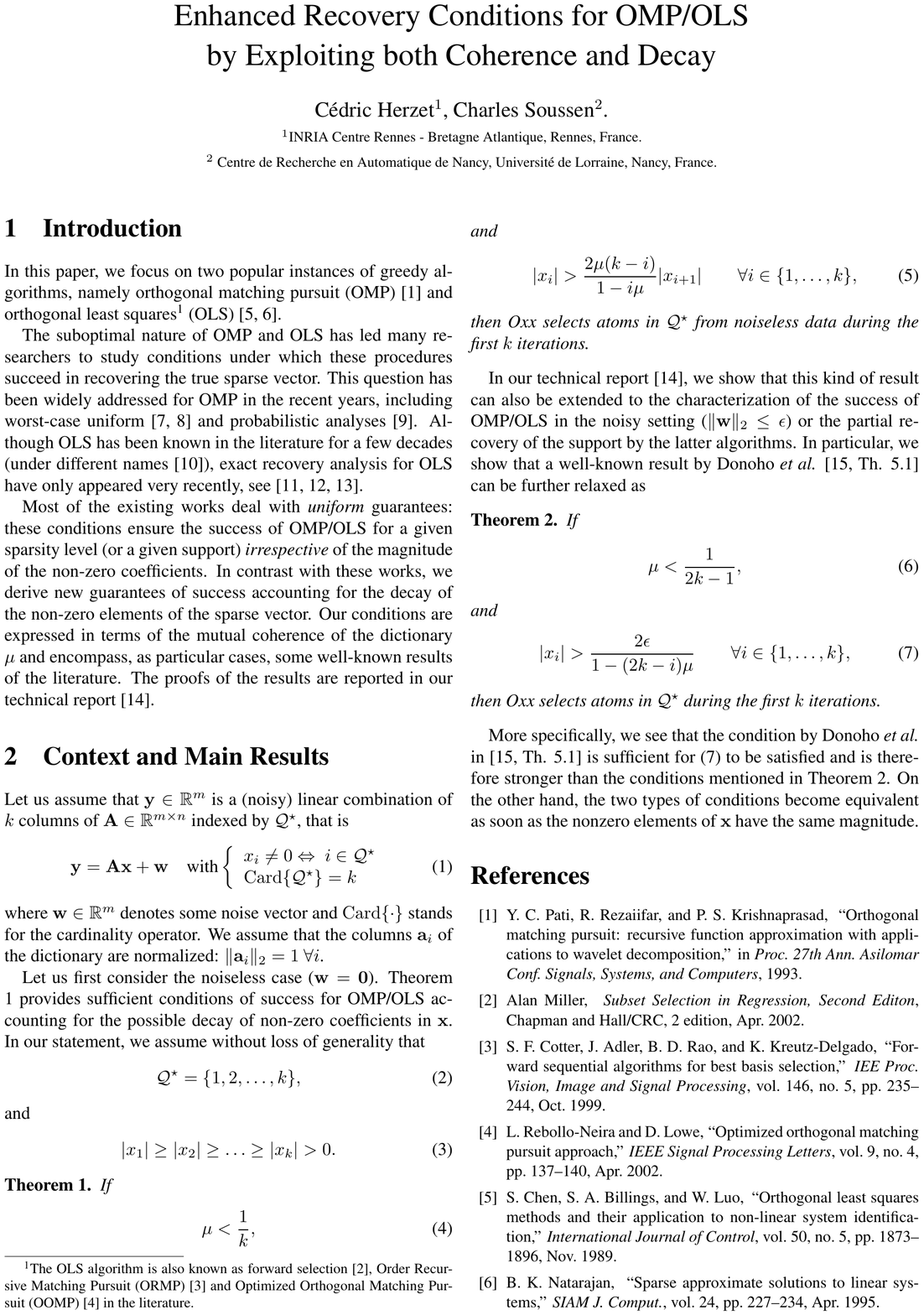}
\label{pdf:Kitic}
\includepdf[offset=1cm -4mm,pages=-,link=true,linkname=Kitic,pagecommand={}]{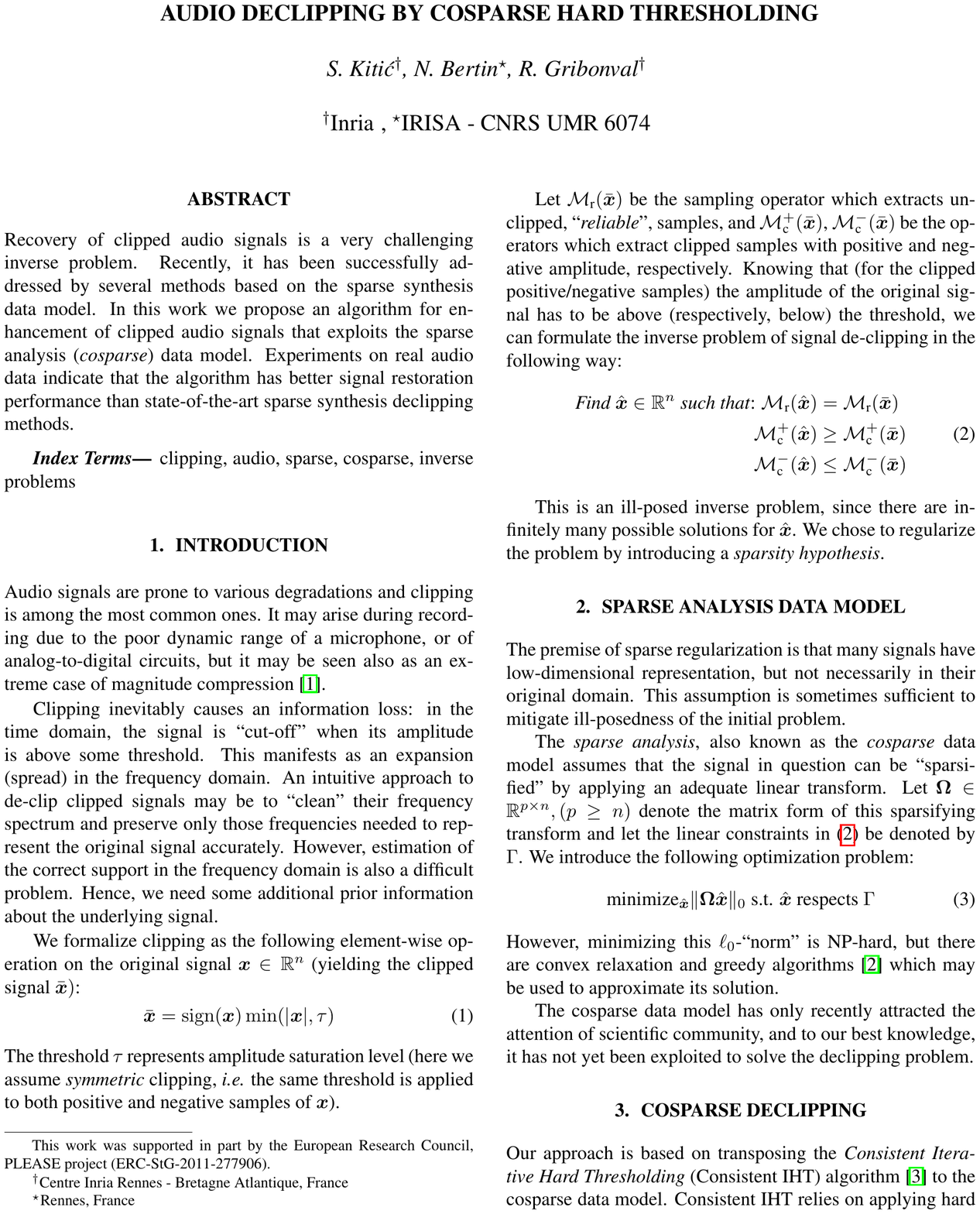}
\label{pdf:LeMagoarou}
\includepdf[offset=1cm -4mm,pages=-,link=true,linkname=LeMagoarou,pagecommand={}]{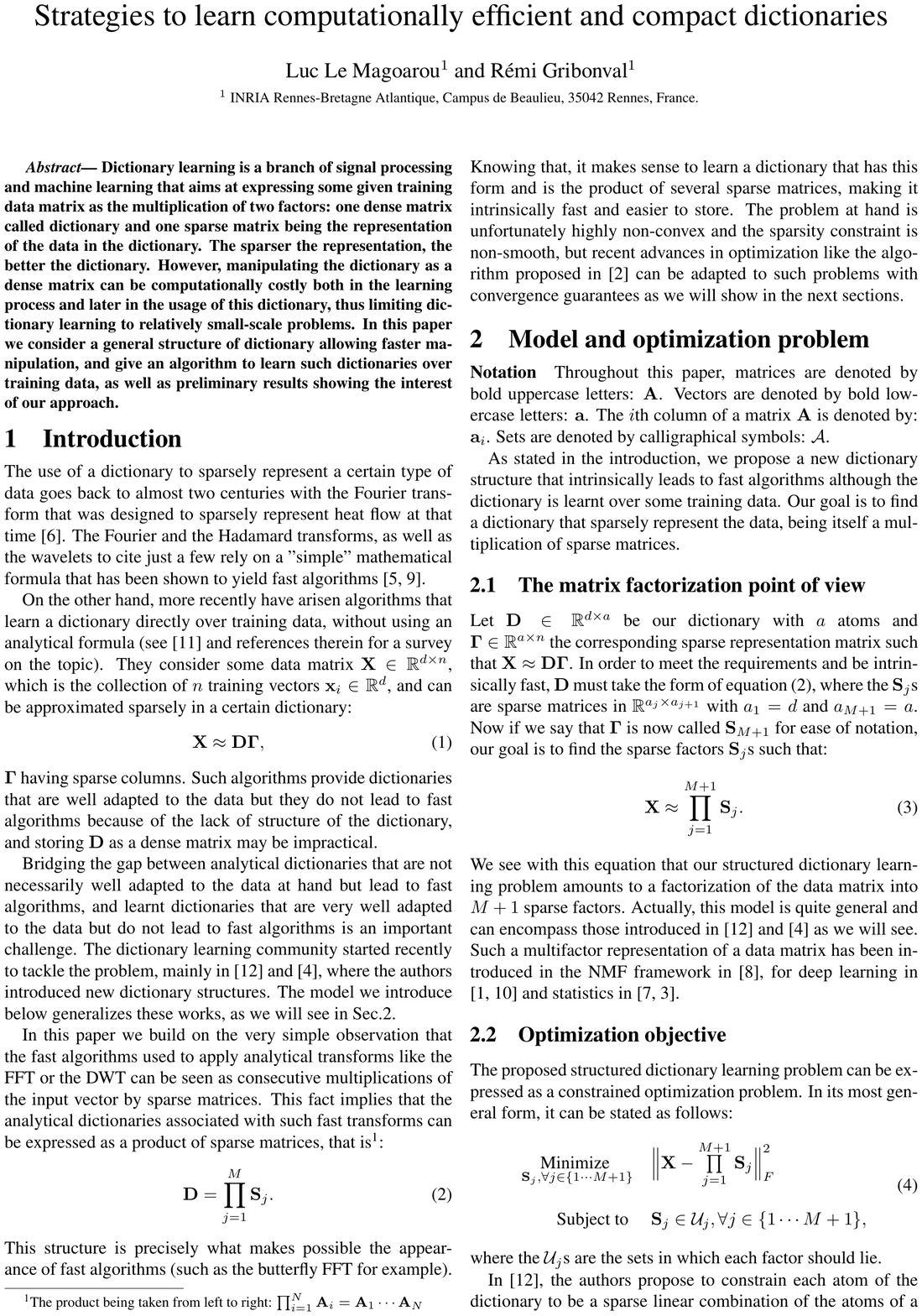}
\label{pdf:Liang}
\includepdf[offset=1cm -4mm,pages=-,link=true,linkname=Liang,pagecommand={}]{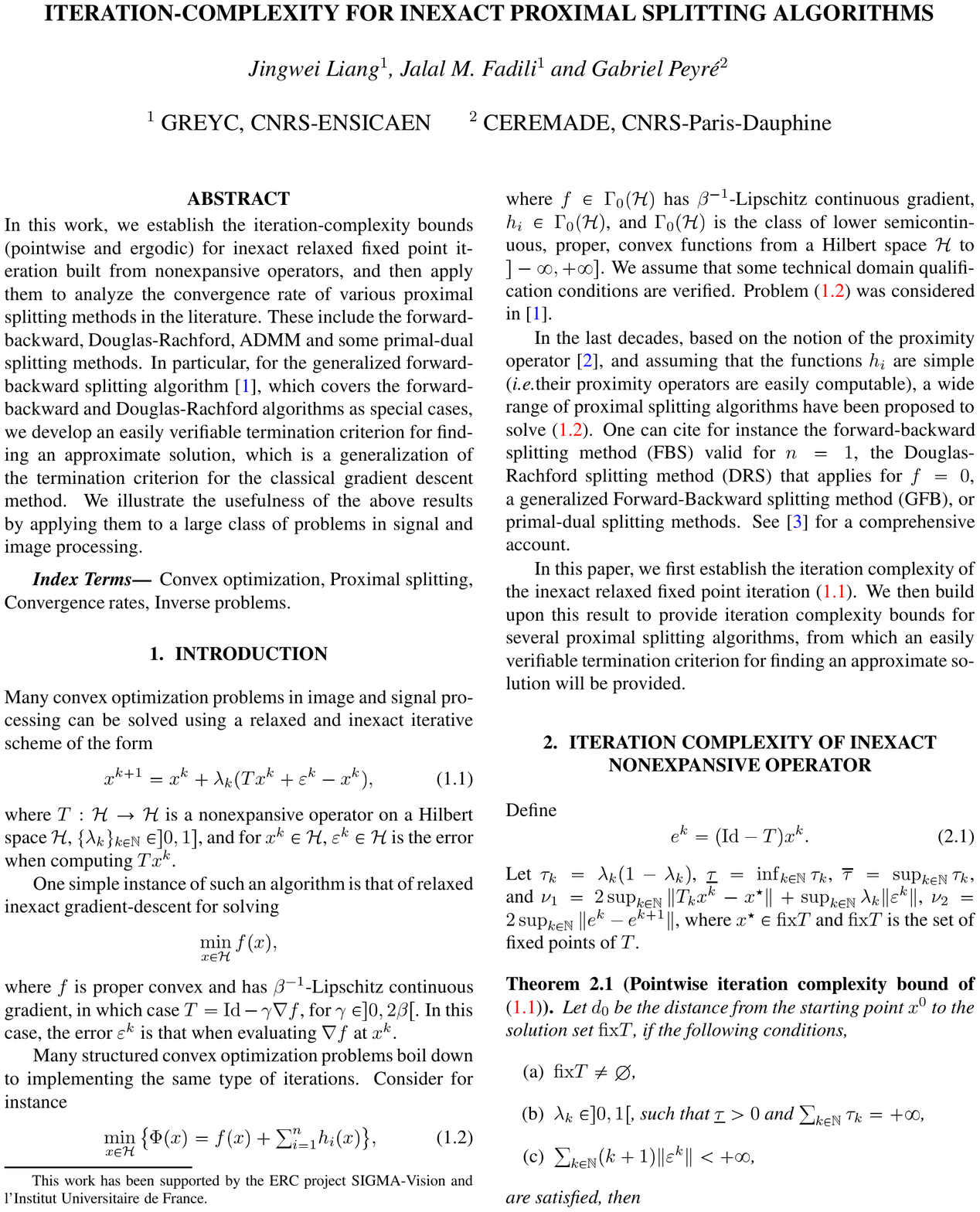}
\label{pdf:Liutkus}
\includepdf[offset=1cm -4mm,pages=-,link=true,linkname=Liutkus,pagecommand={}]{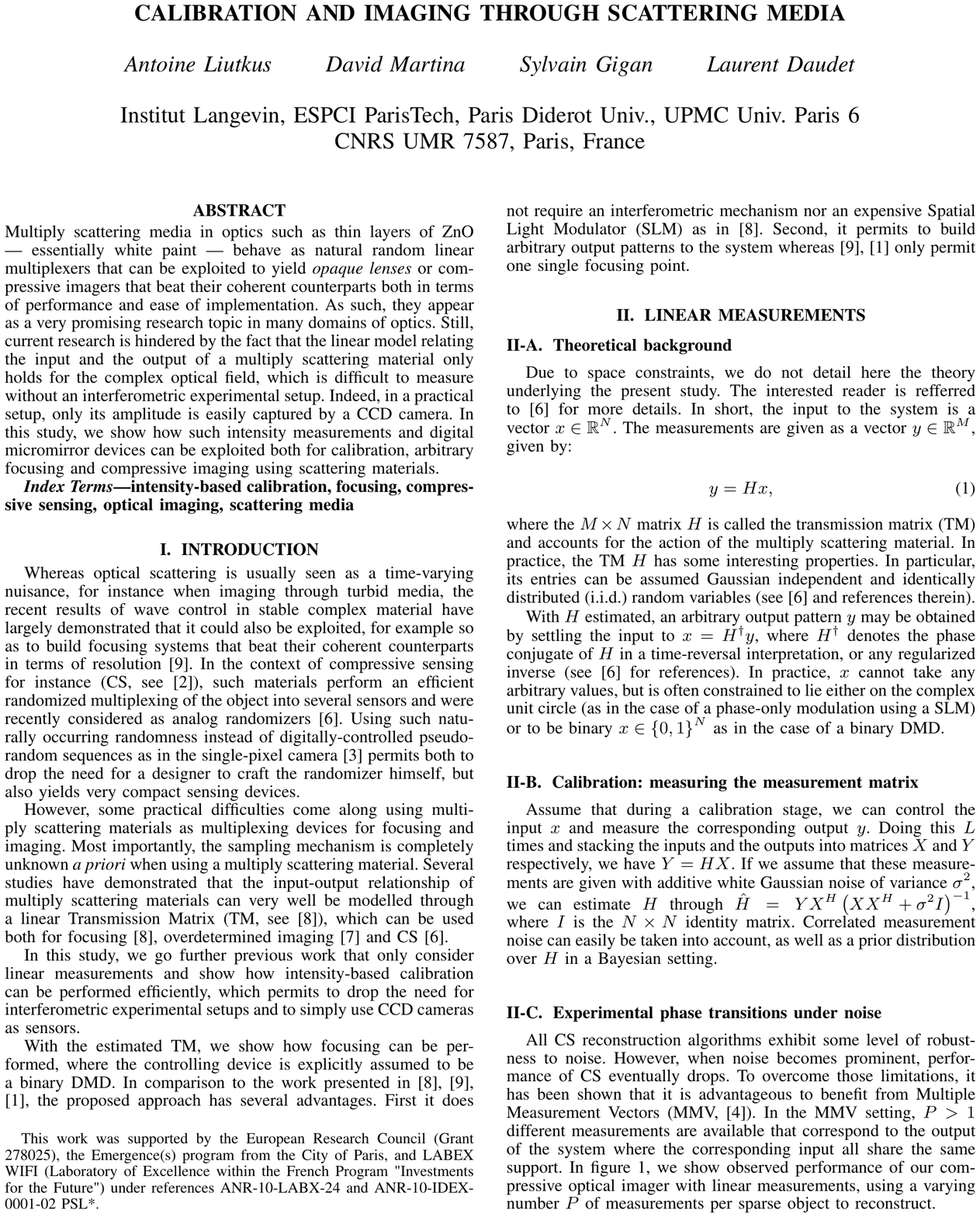}
\label{pdf:Maggioni}
\includepdf[offset=1cm -4mm,pages=-,link=true,linkname=Maggioni,pagecommand={}]{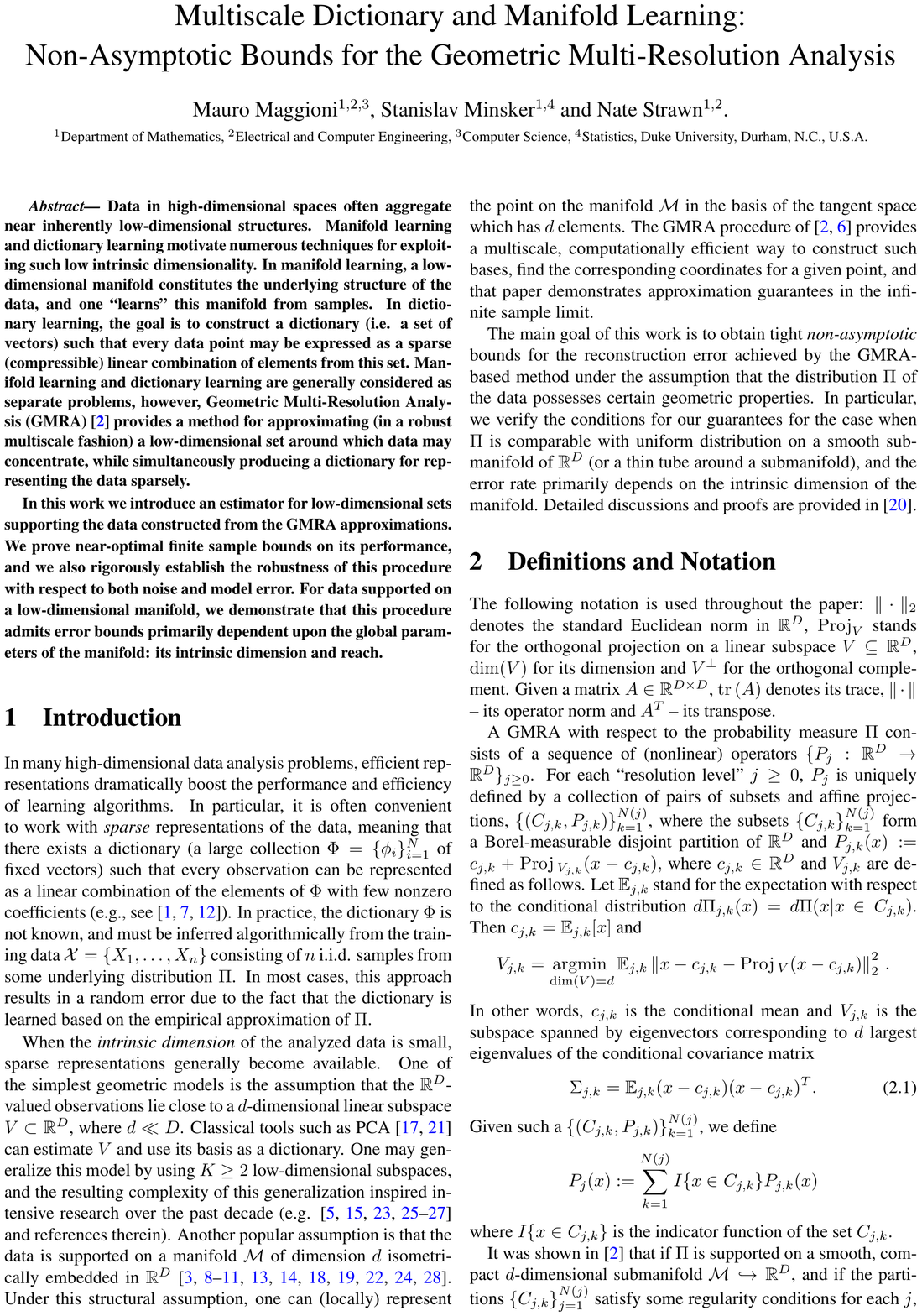}
\label{pdf:Mory}
\includepdf[offset=1cm -4mm,pages=-,link=true,linkname=Mory,pagecommand={}]{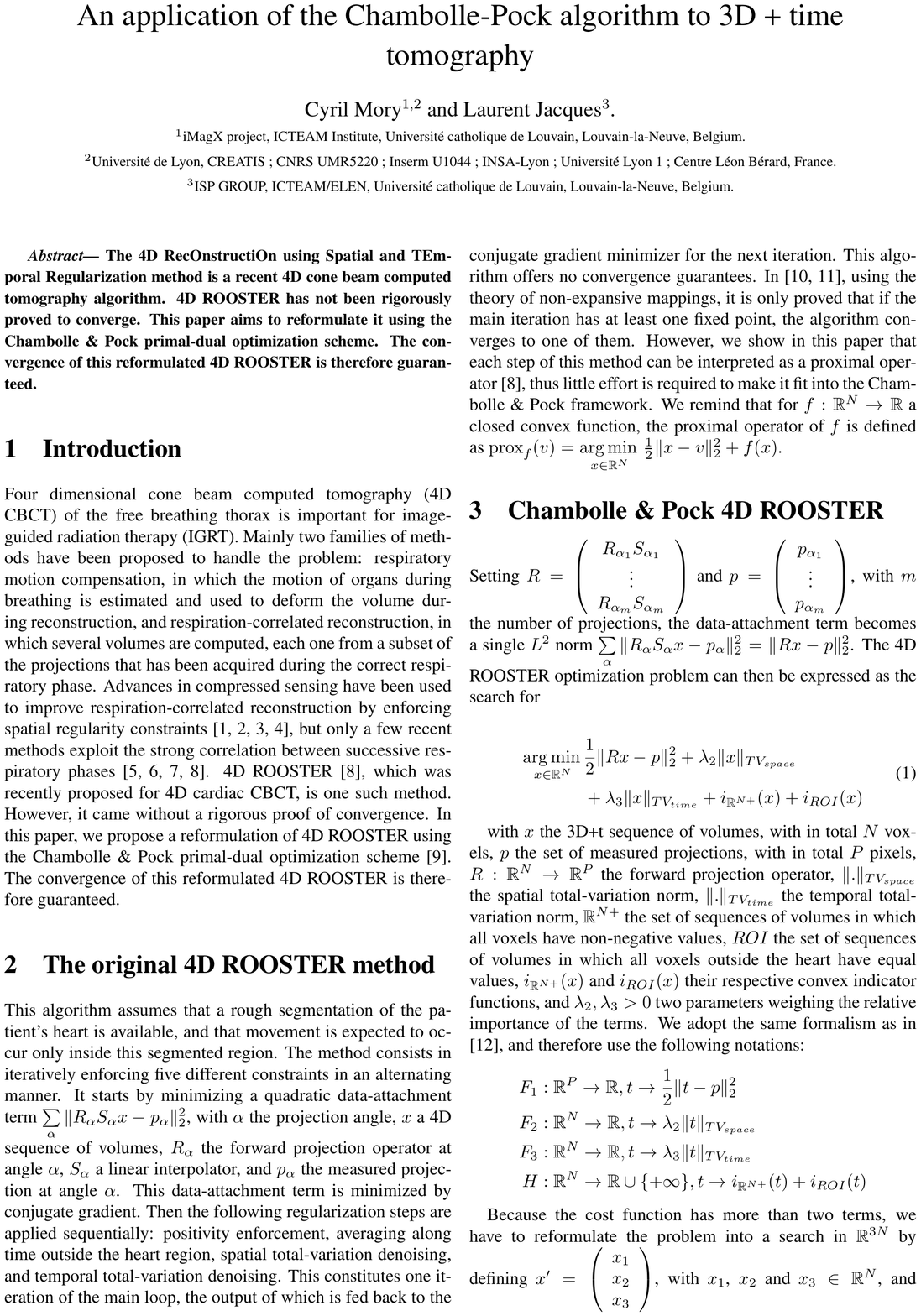}
\label{pdf:Ngole}
\includepdf[offset=1cm -4mm,pages=-,link=true,linkname=Ngole,pagecommand={}]{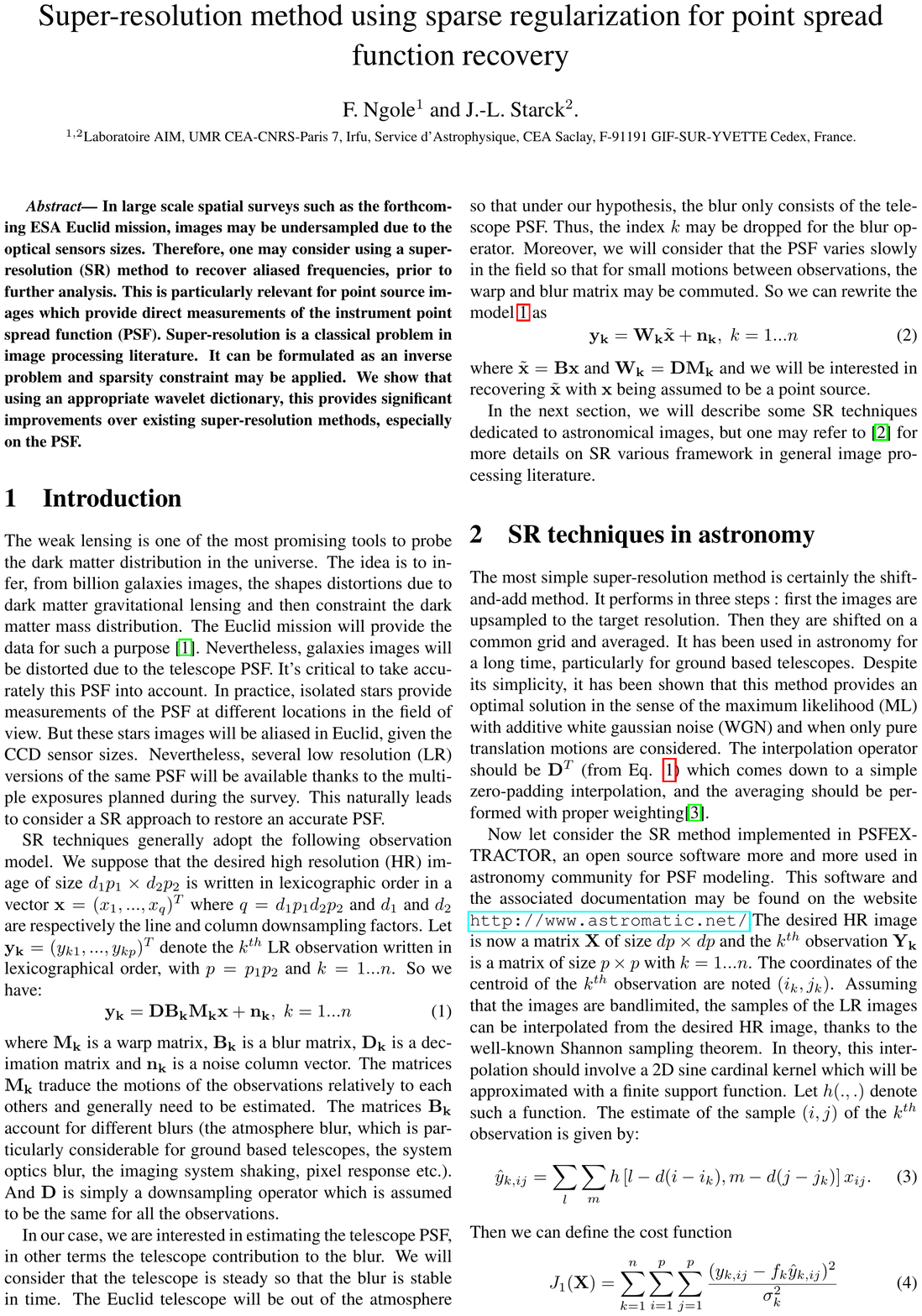}
\label{pdf:Schretter}
\includepdf[offset=1cm -4mm,pages=-,link=true,linkname=Schretter,pagecommand={}]{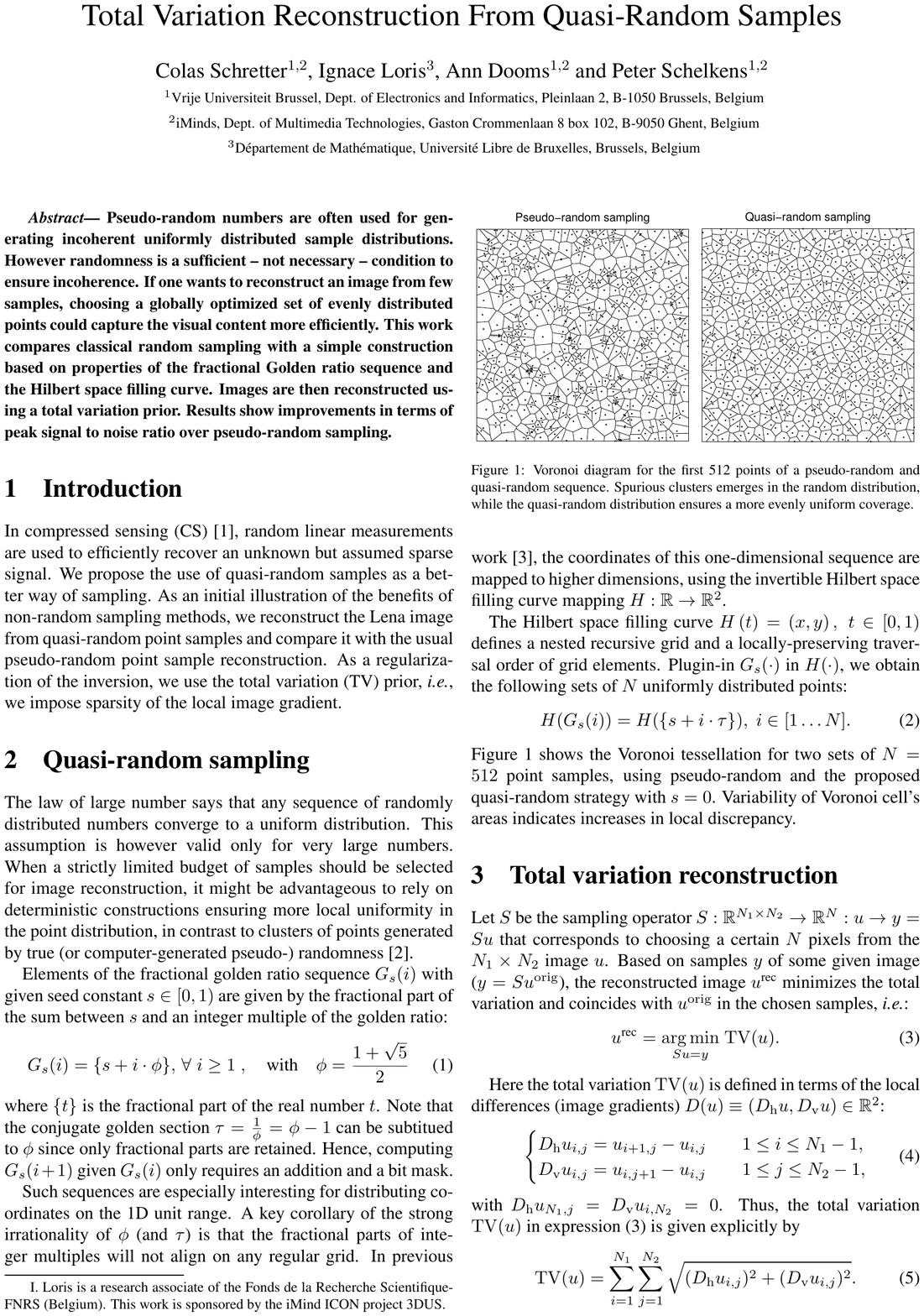}
\label{pdf:Vaiter}
\includepdf[offset=1cm -4mm,pages=-,link=true,linkname=Vaiter,pagecommand={}]{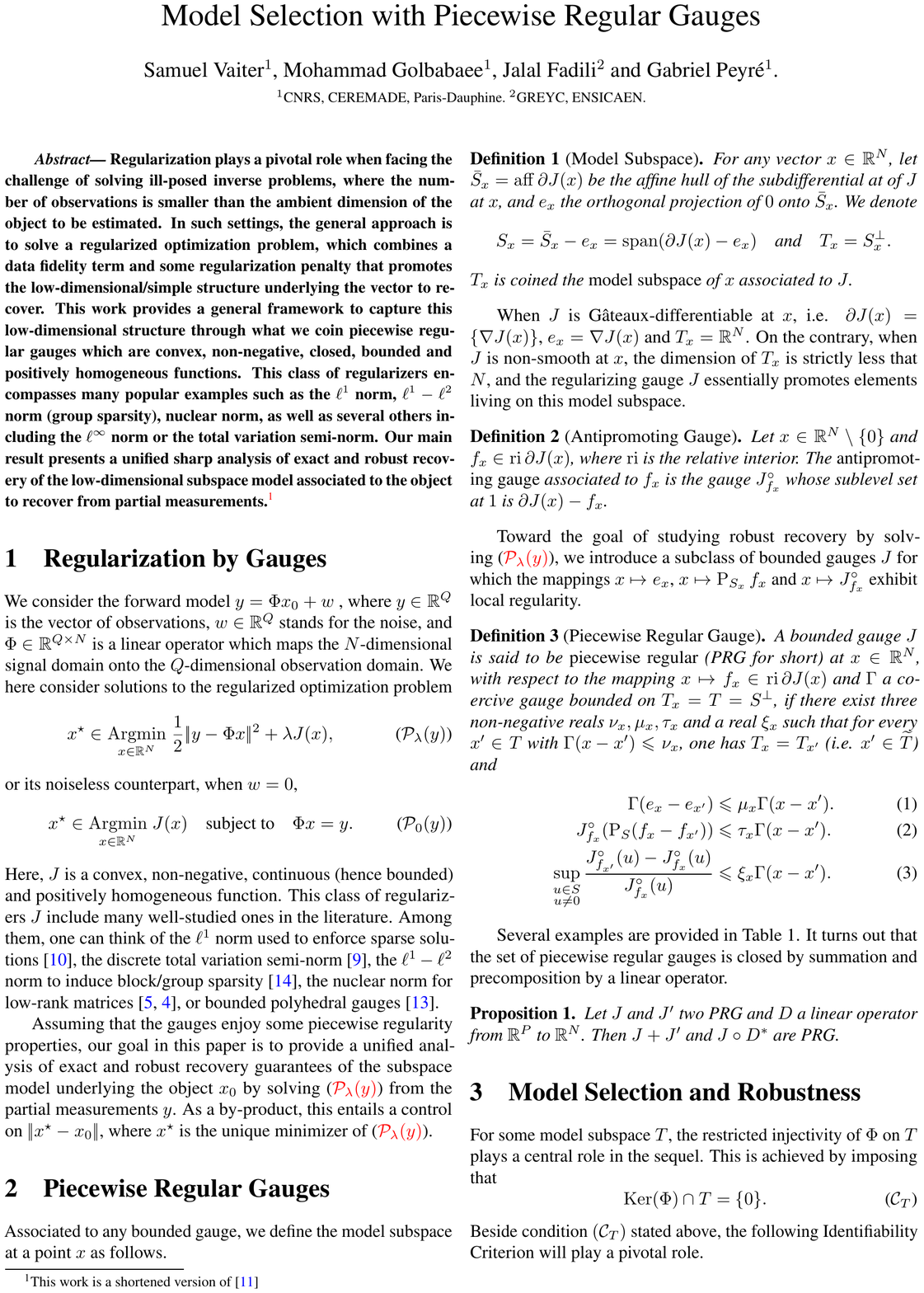}
\label{pdf:Vukobratovic}
\includepdf[offset=1cm -4mm,pages=-,link=true,linkname=Vukobratovic,pagecommand={}]{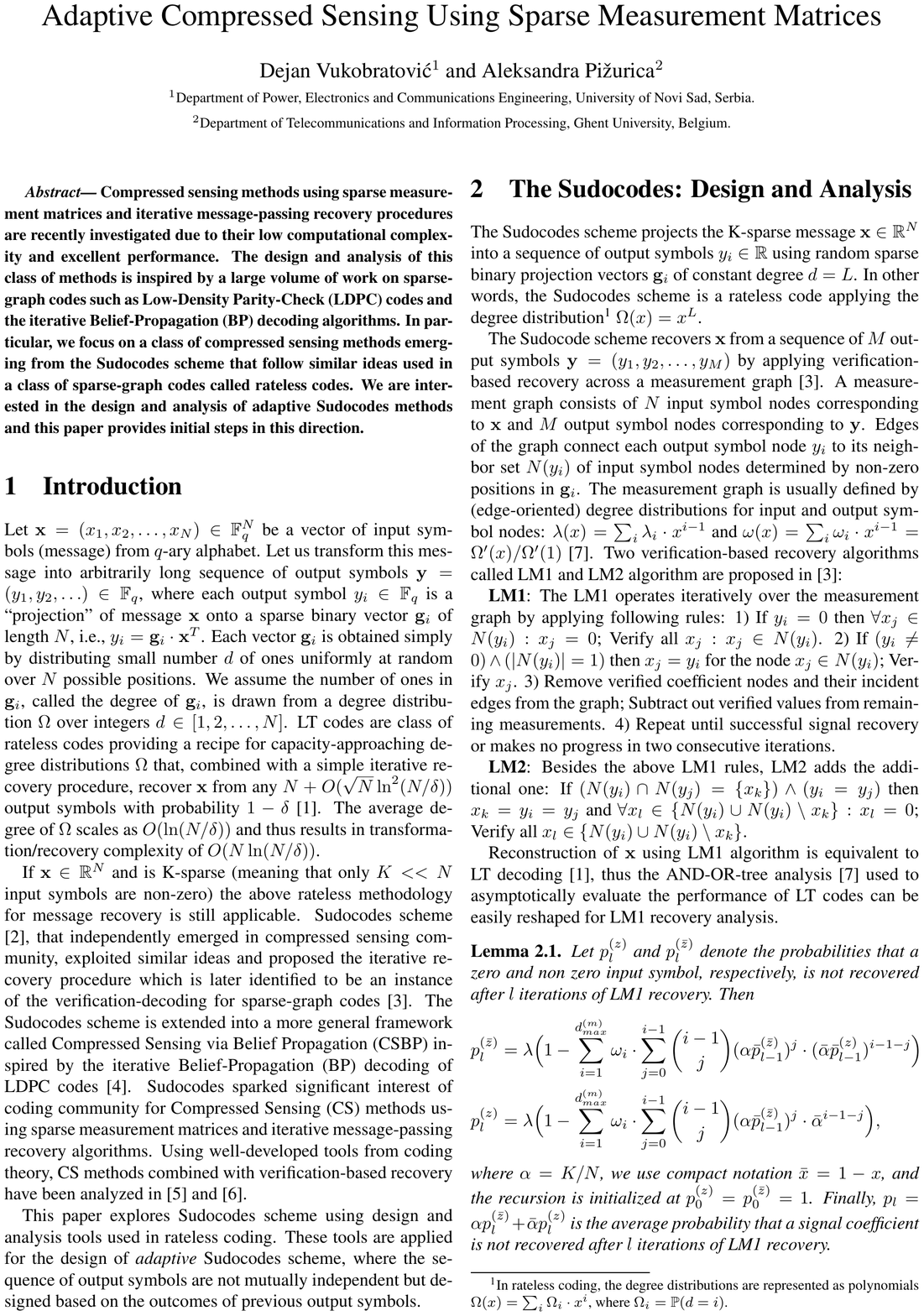}

\newpage
\null
\vfill
\hfill Finished on \today, Louvain-la-Neuve, Belgium.
\null
\end{document}